\documentstyle[11pt,leqno,twoside]{article}
\input xy
\xyoption{all}

\font\es=msbm10 scaled\magstep 1
\def\Rational{\mbox{\es Q}}
\def\Natural{\mbox{\es N}}

\def\proj{\mbox{\es P}}
\def\Proj{\mbox{\es P}}
\def\A{{\cal A}}
\def\dem{{\bf Proof.} \quad}

\begin{document}

\begin{flushleft}
{\Large {\bf LOWER BOUNDS OF THE SLOPE OF FIBRED THREEFOLDS}}\\
\bigskip
\hrule
\end{flushleft}
\vskip 7pt
\begin{flushright}
MIGUEL A. BARJA\

{\scriptsize {
Departament de Matem\`atica Aplicada I.

ETSEIB-Facultat de Matem\`atiques i Estad{\'\i}stica.

Universitat Polit\`ecnica de Catalunya.

Diagonal 64.

08028 Barcelona. Spain.

e-mail: barja@ma1.upc.es.
}}
\end{flushright}

\vglue.2truecm

\begin{flushright}
{\scriptsize {1991 Mathematics Subject
Classification: Primary 14J10, Secondary 14D25}}

\end{flushright}

\vglue.8truecm

\begin{abstract}
We study from a geographical point of view fibrations of threefolds over smooth curves $f: T \longrightarrow
B$ such that the general fibre is of general type.
We prove the non-negativity of certain relative invariants
under general hypotheses and give lower bounds for $K^3_{T/B}$ depending
on other relative invariants. We also study the influence of the relative
irregularity $q(T)-g(B)$ on these bounds. A more detailed study of
the lowest cases of the bounds is given.
\end{abstract}

\baselineskip=14pt

\vglue1truecm

\noindent
{\bf 0. Introduction}

\medskip

We consider fibrations $f:T \longrightarrow B$ from
a complex projective threefold $T$ onto a complex smooth curve $B$. We
always consider $T$ to be normal, with at most canonical singularities
and that $f$ is relatively minimal, i.e. the restriction of $K_{T}$
to any fibre of $f$ is nef.
We are interested in the case where a general fibre $F$ is of general type. Given any fibration
$g:\widetilde T \longrightarrow B$ from a smooth projective threefold
$\widetilde T$  and with fibres of general type
we can always get its relatively minimal associated
fibration by divisorial contractions and flips (see \cite{Ka2},
\cite{M2}).

Our aim is to study $f$ from a geographical point of view and thereby to
relate some numerical invariants
of $T, B$ and $F$. First of all, note that under our assumptions $K_{T}$
(and hence $K_{T/B}=K_{T}-f^{\ast}K_{B}$) is a Weil, $\Rational$-Cartier divisor.
We can consider its associated divisorial sheaves $\omega_{T}$ and
$\omega_{T/B}$, the canonical sheaf of $T$ and the relative ca-

\noindent -----------------

{\footnotesize Partially supported
by CICYT PS93-0790 and HCM project n. ERBCHRXCT-940557}

\newpage

\noindent
nonical
sheaf of $f$, respectively. Then ${\cal E}=f_{\ast}\omega_{T/B}$ is a
locally free sheaf on $B$ of rank $p_{g}(F)$. We then have the well
defined numerical invariants (note that the first one may be a rational
number):

$$K^{3}_{T/B}=K^{3}_{T}-6K^{2}_{F}(b-1)$$
$$\Delta_f=\mbox{deg}{\cal E}$$
$$\chi_f=\chi ({\cal O}_{F})\chi ({\cal O}_{B})-\chi ({\cal O}_{T}).$$

From \cite{F1} we know that ${\cal E}$ is a nef vector bundle and hence
that $\Delta_f \geq 0$. From the nefness of direct images of multiples
of the relative dualizing sheaf, it follows that $K_{T/B}$ is also
nef (see \cite{O}, Theorem 1.4) and so $K^{3}_{T/B} \geq 0$.
Moreover, if $K^{3}_{T/B}=0$, then $f$ is isotrivial (cf. \cite{O}
Corollary 1.5). When $\Delta_f=0$ we can also deduce the isotriviality
of $f$ under some additional hypotheses (see Lemma 1.2).
In general it is not known whether $\chi_f \geq 0$ holds. In \S 1 we
give examples for which $\chi_f$ is negative (see Remark 1.7) and
prove (see Theorem 1.6):

{\bf Theorem 0.1} \quad
{\it
If the Albanese dimension of $T$ is not equal to one, then $\chi _f
\geq 0$ under some extra mild hypotheses.
}

Then we can define for a wide class of fibrations two
different {\it slopes}, $\lambda_1(f)=K^3_{T/B} \Big / \Delta_f
$ (when $\Delta_f \not= 0$) and $\lambda_2(f)=
K_{T/B}^3 \Big / \chi_f$ (when $\chi_f > 0 $) and prove some natural
invariance of them under certain operations.

With these notations we can state the only known general result
on the geography of fibred threefolds over curves, due to Ohno.
A simplified version of the theorem can be stated (see \cite{O} for
a complete reference):

{\bf Theorem 0.2} (\cite{O} Main Theorem 1, Main Theorem 2){\bf .} \quad
{\it Let $f:T \longrightarrow B$ be a relatively minimal fibration
of a threefold over a smooth curve of genus $b$. Assume that a
general fibre $F$ is of general type.
}

(i) {\it Assume $p_g(F) \geq 3$. Then
$K_{T/B}^3 \geq (4 - \varepsilon (p_g(F)))
(\chi ({\cal O}_B)\chi ({\cal O}_F)-
\chi
({\cal O}_T))$, where $\varepsilon (p_g(F))=\frac{4}{p_g(F)}$ or $\frac{8}{
p_g(F)}$ depending on whether $\vert K_F \vert$ is composed with
a pencil or not.}

(ii) {\it If $K_{T/B}^3<4(\chi ({\cal O}_B)\chi ({\cal O}_F) -
\chi({\cal O}_T))$
then $F$ falls in a list of $7$ (possible) families.}

We observe that in \cite{O} it is not proved that $\chi_f$ is non-negative.
In fact, when $\chi_f<0$, Theorem 0.2 gives no information since
$K^3_{T/B} \geq 0$ holds.

In \S 2 we study general lower bounds for $\lambda_1(f)$ and
$\lambda_2(f)$. The main result
(Theorem 2.4) reads:

\newpage

{\bf Theorem 0.3.} \quad
{\it With the same hypotheses of Theorem 0.1, if $p_g(F)\geq 3$ and
$\chi_f >0$, then $\lambda_2(f)\geq (9-\widetilde {\varepsilon}(p_g(F)))$ except if
$F$ is fibred by hyperelliptic, trigonal or tetragonal curves, or
$|K_F|$ is composed, where ${\widetilde \varepsilon}(p_g(F)) \sim O(\frac
{1}{p_g(F)})$.}

In fact Theorem 2.4 is much more explicit and gives extra bounds for
the exceptional cases.

Section \S 3 is devoted to the study of the influence of the
irregularity of $T$ on the slope.
In the study of fibred surfaces we have the general inequality
$\lambda(f)\geq 4 -
\varepsilon (F)$ and that $\lambda(f)\geq 4$ when $q(S)>g(B)$
(cf. \cite{X1}).
Similarly we get for threefolds
(cf. Theorem 3.3):

{\bf Theorem 0.4.} \quad {\it If $q(T)>g(B)$, then
$\lambda_2(f)\geq 9$ except when $F$ has an irrational pencil of
hyperelliptic, trigonal or tetragonal curves.}

Theorem 3.3 also gives explicit bounds and a structure result in the
exceptional cases. The key point here is to use the condition
$q(T) > g(B)$ to construct new fibrations with the same slope and
with fibres of higher invariants. The result follows then from a limit
process in Theorem 0.3.
Also we remark that, following \cite{BT}, the hypothesis
$q(T)>g(B)$ can be weakened to the condition that ${\cal E}=f_{\ast}\omega_{T/B}$ has
a locally free, rank one, degree zero quotient.

Finally in \S 4 we study fibrations with very low slope ($\lambda_{2}<4$).
These are known to exist (cf. \cite{O} p.664);
in \cite{O}, Ohno gives a classification
of them in seven possible families as stated in Theorem 0.2 (see
Theorem 4.1 below for a
complete description). We prove (see Theorem 4.2):

{\bf Theorem 0.5.} \quad
{\it
Let $f:T\longrightarrow B$ be a
relatively minimal fibration of a normal,
projective threefold $T$ with only canonical singularities onto a smooth
curve $B$ of genus $b$. Assume that a general fibre $F$ is of general
type with $p_{g}(F)\geq 3$ and that $\chi
_{f}=\chi({\cal O}_{F})\chi({\cal O}_{B})-\chi({\cal O}_{T})>0$.

Then, if  $\lambda _{2}(f)<4$, we have:
\begin{enumerate}
\item[{\rm (i)}]  $q(T)=b$
\item[{\rm (ii)}]  ${\cal E} =f_{\ast }\omega_{T/B}$  has no invertible degree
zero quotient sheaf (in particular,  ${\cal E} $ is ample provided
$b\leq 1$).
\item[{\rm (iii)}] If $p_g(F) \geq 15$, then $F$ has a rational pencil of curves
of genus $2$.

\item[{\rm (iv)}] If $p_g(F) \leq 14$, then one of the following holds:
\begin{itemize}
\item[{\rm (a)}]  $F$ has a rational pencil of hyperelliptic curves.
\item[{\rm (b)}]  $F$ has a rational pencil of trigonal curves and $q(F)=0$.
\item[{\rm (c)}]  $F$ is a quintic surface in $\Proj ^{3}$.
\end{itemize}
\end{enumerate}
}

\newpage

In fact we have a more concrete description of case (iv) (b). Case
(iv) (c) is doubtful to exist (see Remark 4.3).

\vglue.5truecm

The author wishes very much to thank Professor Kazuhiro Konno who pointed
out the problem to him during his visit to Barcelona and for his
encouragement and advice. He also wants to thank his advisor, Professor
Juan Carlos Naranjo, for help and continuous support.

\vglue.5truecm

{\bf Notations and conventions}

We work over the field of complex numbers. Varieties are always
assumed to be projective, reduced and irreducible. Symbols $\sim$,
$\sim _{\Rational}$ and $\equiv$ stand for linear equivalence,
$\Rational$-linear equivalence and numerical equivalence respectively.

If ${\cal F}$ is a coherent sheaf on a variety $X$ we usually put
$\chi ({\cal F})$ instead of $\chi (X,{\cal F})$.

\vglue.5truecm

\noindent
{\bf  1. Slopes of fibred threefolds}

\medskip

{\bf Definition 1.1.} \quad
Let $f:T \longrightarrow B$ be a fibration of a normal,
projective threefold with only canonical
singularities onto
a smooth curve. Let $F$ be a general fibre of $f$ and put
$b=g(B)$. We define
\begin{eqnarray*}
&&\Delta_f=\mbox{\rm deg }f_{\ast }w_{T/B}\\
&&\chi_f =\chi ({\cal O}_{B})\chi ({\cal O}_{F})-\chi ({\cal O}_{T})
\end{eqnarray*}

{\bf Lemma 1.2.} \quad
{\it
\begin{enumerate}
\item[{\rm (i)}]  $\Delta_f=\chi_f+\mbox{\rm deg}R^1f_{\ast}\omega_{T/B}\geq
\chi_f$
\item[{\rm (ii)}]  $\Delta_f\geq 0$. If $\Delta_f=0$ and $|K_{F}|$ is
birational, then  $f$ is isotrivial.
\item[{\rm (iii)}] If $\beta:\widetilde T \longrightarrow T$ is a nonsingular
model of $T$ and $\widetilde f=f \circ \beta$, then
$\chi_f=\chi_{\widetilde f}$, $\Delta_f = \Delta_{\widetilde f}$.
\end{enumerate}
}

\dem

\begin{enumerate}
\item[(i)] Follows from \cite{O} Lemma 2.4 and 2.5.
\item[(ii)]  $\Delta_f\geq 0$ follows from the nefness of ${\cal E} $
(\cite{F1}). If $\Delta_f=0$ and $|K_{F}|$ is birational we can apply
\cite{Ko4} I (see also \cite{M}, 7.64).
\item[(iii)] Canonical singularities are rational (cf. \cite{El})
and hence $R^i\beta_{\ast}{\cal O}_{\widetilde T}=0$ for $i \geq 1$
(cf. \cite{Ke}, p.50).
Hence $\chi ({\cal O}_{\widetilde T})=\chi ({\cal O}_T)$. The same
holds for general fibres $F$ and $\widetilde F$ of $f$ and $\widetilde f$
respectively, so $\chi_f=\chi_{\widetilde f}$.

By Grauert-Riemenschneider's vanishing we have
$R^i\beta_{\ast}\omega_{\widetilde T}=0$
for $i \geq 1$. Hence using the spectral sequence $E_2^{p,q}=R^pf_{\ast}
(R^q\beta_{\ast}\omega_{\widetilde T}) \Rightarrow R^{p+q}\widetilde
f_{\ast}\omega_{\widetilde T}$ we obtain that for every $i \geq 0$,
$R^if_{\ast}\omega_T=R^i \widetilde f_{\ast}\omega_{\widetilde T}$ holds.
\end{enumerate}
\quad $\Box$

{\bf Definition 1.3.} \quad
With the above notations if $f$ is relatively minimal
and a general fibre $F$ of $f$ is of general type we define:
\begin{eqnarray*}
&&\lambda
_{1}(f)=\displaystyle\frac{K^{3}_{T/B}}{\Delta_f}\quad\mbox{if}\quad
\Delta_f\not= 0 \\
&&\lambda
_{2}(f)=\displaystyle\frac{K^{3}_{T/B}}{\chi_f}\quad\mbox{if}\quad
\chi_f> 0
\end{eqnarray*}

{\bf Remark 1.4.} \quad
In the case of fibrations of surfaces over
curves we actually have  $\Delta_f=\chi ({\cal O}_S) -\chi ({\cal O}_B)
\chi ({\cal O}_F)=-\chi_f\geq 0$ (the minus sign here is a matter
of the dimension of the variety) and vanishing
holds in
the locally trivial case. Then when  $f$ is not locally trivial we define
$\lambda (f)=K^{2}_{S/B}\Big / \Delta_f$.

Here we have two different possibilities for the slope of
$f:K^{3}_{T/B}\Big /\Delta_f$ or $K^{3}_{T/B}\Big /\chi_f$.
As we will see in \S 2, natural methods provide lower bounds for
$\lambda_1(f)$ (hence also for $\lambda_2(f)$: see Lemma 1.5 below). Note
that from the {\it geographical} point of view the most interesting one
is $\lambda_2(f)$. But for this choice we do not know whether
$\chi_f \geq 0$. The aim of this section is to show that this
actually happens for general fibrations.

{\bf Lemma 1.5.} \quad
{\it
Assume  $\chi_f>0$. Then
\begin{enumerate}
\item[{\rm (i)}]  $\lambda _{2}(f)\geq \lambda _{1}(f)$
\item[{\rm (ii)}] If $\sigma :\widetilde{B}\longrightarrow B$  is an
\'etale map
and
$\widetilde{f}:\widetilde{T}=T\displaystyle\mathop{\times
}_{B}\widetilde{B}\longrightarrow \widetilde{B}$ is the induced
fibration, then $\lambda _{i}(f)=\lambda _{i}(\widetilde{f})$,  $i=1,2$.
\item[{\rm (iii)}] If $\widetilde{T}\displaystyle\mathop{\longrightarrow
}^{\alpha }T$ is an \'etale map such that  $\widetilde{f}=f\circ\alpha $
has connected fibres, then $\lambda _{2}(f)=\lambda _{2}(\widetilde{f})$
(but in general $\lambda _{1}(f)\not= \lambda _{1}(\widetilde{f})$).
\end{enumerate}
}

\dem
(i) is obvious by Lemma 1.2 (i).

(ii) If $\sigma $ does not ramify over the images of singular fibres of
$f$ then $\widetilde{T}=T\displaystyle\mathop{\times }_{B}\widetilde{B}$
is again a normal, relatively minimal threefold over
$\widetilde B$ with only canonical singularities (cf. \cite{M}, 4.10).
Clearly $K_{\widetilde T/\widetilde B}^3 =(\mbox{deg}\sigma)K_{T/B}^3$
and deg$R^i\widetilde f_{\ast} \omega_{\widetilde T/\widetilde B}=
(\mbox{deg}\sigma)R^if_{\ast}\omega_{T/B}$ by flat base change. Then
$n\Delta_f=\Delta_{\widetilde f}$, $n\chi_f=\chi_{\widetilde f}$ and
we are done.

(iii)
By \cite{Fu} Ex. 18.3.9 we have $K_{\widetilde T}=\alpha^{\ast}(K_{T})$
and that $\chi({\cal O}_{\widetilde T})=({\rm deg} \alpha)
\chi ({\cal O}_{T})$. \quad $\Box$

The question now is whether  $\chi_f\geq 0$ holds. This is not always
true (see Remark 1.7). We give criteria for its non-negativity depending
on the Albanese dimension of $T$. By Lemma 1.2 we have
$\chi_f=\chi_{{\widetilde f}}$ where ${\widetilde f}=f \circ \beta$,
$\beta : {\widetilde T} \longrightarrow T$ being a desingularization.
Hence we can assume $T$ is smooth.

First of all, consider $t\in B$, such that  $F_{t}$ is smooth, and the
Albanese maps
$$
\xymatrix{
F_{t}\ar[r]^{\mbox{alb}_{F_{t}}}\ar[d]^{i_{t}}&\mbox{Alb}(F_{t})\ar[d]
_{(i_{t})_{\ast }}\\
T
\ar[r]^{\mbox{alb}_{T}}\ar[d]^{f}&\mbox{Alb}(T)\ar[d]
_{f_{\ast }}\\
B
\ar[r]^{\mbox{alb}_{B}}&\mbox{Alb}(B)}
$$
and let  $\Sigma =\mbox{alb}_{T}(T)$. We set  $a={\rm dim} (alb_{T}(T))=
\dim\Sigma $. Let $S\longrightarrow \Sigma $ be a minimal
desingularization of $\Sigma $ and $\pi :\widetilde{T}\longrightarrow S$
the induced map on a birational model of $T$.

Note that by rigidity,  $\mbox{Im }(i_{t})_{\ast }=A$ is an
abelian variety independent of $t$, of dimension $q(T)-b$.

Also consider the induced map
${\rm Pic}^0T\displaystyle\mathop{\longrightarrow
}^{(i_{t})^{\ast }} {\rm Pic}^0(F_{t})$, whose image is
$\widehat{A}\hookrightarrow  {\rm Pic}^0(F_{t})$. We say that  $f$
is {\it special} if for general  $t\in B$, $\widehat A
\hookrightarrow (h_{t})^{\ast
}({\rm  Pic}^0C_{t})\subseteq {\rm Pic}^0(F_{t})$, for some
$h_{t}:F_{t}\longrightarrow C_{t}$ a fibration over a curve of genus
$g(C_{t})\geq 2$. {\cal O}therwise we say that $f$ is {\it general}.

{\bf Theorem 1.6.} \quad
{\it
Let $f:T\longrightarrow B$ be a
fibration of a normal, projective threefold
with only canonical singularities onto a smooth curve of genus $b$. Let
$F$ be a general fibre of $f$. Let $a=\mbox{\rm dim}\,alb_{T}(T)$.

Then  $\chi_f=\chi ({\cal O}_{B})\chi ({\cal O}_{F})-\chi ({\cal O}_{T})\geq 0$ provided one of the
following conditions holds:
\begin{enumerate}
\item[{\rm (i)}]  $b \leq 1$ and $\chi ({\cal O}_{T}) \leq 0$.
\item[{\rm (ii)}]  $a=2$, $b \geq 1$ and $h^{\circ }(S,\pi_{\ast }w_{\widetilde T/S})\not= 0$.
\item[{\rm (iii)}] $a=3$,  $f$ is special and is semistable.
\item[{\rm (iv)}] $a=3$ and  $f$ is general.
\end{enumerate}
}

{\bf Remark 1.7.} \quad
Part (i) of the theorem is trivial since when $b \leq 1$, $\chi_f \geq
- \chi {\cal O}_T$.
We only want to remark that condition $\chi ({\cal O}_T) \leq 0$ holds in
most cases. Indeed, if $T$ is smooth and $K_T$ ample, then
$\chi ({\cal O}_T)=\frac{1}{24}c_1c_2<0$ by Miyaoka-Yau inequality.
Also, if $T$ is minimal and Gorenstein, $\chi ({\cal O}_T) \leq 0$ holds
(cf. \cite{M1}). Finally, if $a=3$, then $\chi ({\cal O}_T)\leq 0$ by a
consequence of generic vanishing results (see \cite{GL1}). Observe
that if $a=0$, then necessarily $q(T)=b=0$ and hence this possibility
is included in (i), i.e., we need to know whether $\chi ({\cal O}_T)
\leq 0$.
Extra conditions included for the cases $a=2,3$ are very mild. This
is clear for the case $a=3$. In the case $a=2$ observe that if $E$
is a general curve on $S$ and $H$ is its pullback on ${\widetilde T}$,
then $h^0(S,\pi_{\ast}\omega _{{\widetilde T}/S})=h^0(E, \pi_{\ast}
\omega _{H/E})$ (see proof of part (iii) of the theorem). Then
by Proposition 1.8.(i) $\pi _{\ast} \omega_{H/E}$ is nef
and (in general) contains an ample vector bundle, so condition $h^0(E,\pi_{\ast}
\omega_{H/E}) \neq 0$ is not very restrictive.

We want to stress that in the statement of the theorem some hypotheses
are needed since $\chi_f \geq 0$ {\it is not always true}. Indeed, we
can construct counterexamples following \cite{M1} Remark 8.7. Let $(C_i, \tau_i)$
(i=1,2,3) be smooth
curves with an involution. Let $D_i=C_i/\tau_i$ and $g_i=g(C_i) \geq 2$,
$b_i=g(D_i)$. Consider $X=C_1 \times C_2 \times C_3$ and
$\tau : X \longrightarrow X$ the involution acting on $C_i$ as $\tau_i$.
Consider $T=X/\tau$. Then $T$ is a threefold of general type with a
finite number of (canonical) singularities, endowed with a fibration
$f: T \longrightarrow D_1=:B$ with general fibre $F \cong C_2 \times
C_3$ (hence it is isotrivial). Then:

$$\begin{array}{rl}
\chi ({\cal O}_B)&=1-b_1\\
&\\
\chi ({\cal O}_F)&=(1-g_1)(1-g_2)\\
&\\
h^{1}(T,{\cal O}_T)&=q(T)=b_1+b_2+b_3\\
&\\
h^{2}(T,{\cal O}_T)&=b_1b_2+b_1b_3+b_2b_3+(g_1-b_1)(g_2-b_2)+
(g_1-b_1)(g_3-b_3)+\\
&\\
&\quad +(g_2-b_2)(g_3-b_3)\\
&\\
h^{3}(T,{\cal O}_T)&=b_1b_2b_3+(g_1-b_1)(g_2-b_2)b_3+(g_1-b_1)(g_3-b_3)b_2+\\
&\\
&\quad +(g_2-b_2)(g_3-b_3)b_1.
\end{array}$$

If we take $b_1=b_2=b_3=0$ ($C_i$ must be hyperelliptic then) we
obtain $a=0$ and $\chi_f <0$. Any base change, \'etale over the critical
points of this fibration, to
a curve of positive genus produces a new fibration with
$q({\widetilde T})={\widetilde b} \geq 1$ (hence with $a=1$) and
$\chi_{\widetilde f} <0$.

If we take $b_1=b_2=1$, $b_3=0$ we obtain again $\chi_f <0$ and
$q(T)>b$ (so $a \geq 2$).

Finally we must say that we do not have any
reasonable criteria for the nonnegativity of $\chi
_{f} $ when $a=1$, which corresponds to the case when $q(T)=b$,
$f_{\ast}=\mbox{\rm Id }$ and  $\mbox{\rm alb}_{T}$ factors through $f$
($f$ is, then, an Albanese fibration).
Nevertheless this is precisely the case in which we are not interested
in Theorem 3.3.

\dem
(ii) From Remark 2.3 we have that $\chi_f=\chi_{\widetilde f}$, where
$\widetilde{f}=f\circ \beta $, $\beta :\widetilde{T}\longrightarrow T$
is a desingularization.
Hence we can assume  $T$ smooth. We can always assume, by the
same arguments, that $\pi :T\longrightarrow S$ has branch locus
contained in a normal crossings divisor.

Since $b \geq 1$ we now have a factorization of $f$
$$
\xymatrix{
T\ar[dd]_{f}\ar[dr]^{\pi }&{}\\{}&S\ar[dl]^{g}\\ B&{}
}
$$
where  $g$ need not to be a relatively minimal fibration. Let
$C_{t}=g^{-1}(t)$ be a general fibre and  $\pi
_{t}:F_{t}\longrightarrow C_{t}$ the induced fibration. In order to
simplify the notation we sometimes will use $C$ instead of
$C_{t}$.
Let $G$ be a
general fibre of $\pi _{t}$. Note that we have $(\pi _{t})_{\ast}
\omega_{F_{t}/C_{t}}=\pi _{\ast }\omega_{T/S}\otimes {\cal O}_{C_{t}}$
and hence $(R^{1}\pi
_{t} )_{\ast } {\cal O}_{F_{t}}=(R^{1}\pi _{\ast }{\cal O}_{T})\otimes
{\cal O}_{C_{t}}$ by relative duality on $C_t$.

Take a  $n$-torsion element  ${\cal L}\in  {\rm Pic}^0(S)$ verifying that for
$1\leq i\leq n-1$ \break  ${{\cal L}}^{\otimes i}\not\in
g^{\ast }( {\rm Pic}^0(B))$ (this is possible since $S$ is of Albanese
general type by construction) and such that $h^{\circ}(C_{t},
R^{1}(\pi _{t})_{\ast }{\cal O}_{F_{t}}\otimes {\cal L}_{|C_{t}})=0$ (this is also
possible since $\{\widetilde{{\cal L}}\in  {\rm Pic}^0(C_{t})\, | \, h^{\circ
}(C_{t},(R^{1}\pi _{t})_{\ast }{\cal O}_{F_{t}}\otimes\widetilde{{\cal L}})\not= 0\}$
is a finite set (see Proposition 1.8) and the image of
$ {\rm Pic}^0(S)\longrightarrow  {\rm Pic}^0(C_{t})$ is a subtorus of
positive dimension otherwise $q(S)=b$, a contradiction).

Let ${\cal M}=\pi ^{\ast }{\cal L}\in  {\rm Pic}^0(T)$. Since $\pi$ has a normal crossings
ramification locus, we have that $R^1\pi_{\ast}{\cal O}_T$ is locally free
(cf. \cite{Ko2} Theorem 2.6 and \S 3) and hence $g_{\ast}(R^1\pi_{\ast}{\cal O}_T
\otimes {\cal L})$ is torsion free (hence it is locally free since $B$
is a smooth curve). Then:
$$g_{\ast }(R^{1}\pi _{\ast }{\cal O}_{T}\otimes
{\cal L})=0$$
since $\mbox{rk }g_{\ast }(R^{1}\pi _{\ast
}{\cal O}_{T}\otimes{\cal L})=h^{\circ}(C_{t},R^{1}(\pi _{t})_{\ast
}{\cal O}_{F_{t}}\otimes{\cal L}_{|C_{t}})=0$ by the choice of ${\cal L}$.

Using the spectral sequence $E^{p,q}_{2}=R^{p}g_{\ast }(R^{q}\pi
_{\ast }{\cal F})\Rightarrow R^{p+q}f_{\ast }{\cal F}$ and that $R^{2}\pi
_{\ast }{\cal O}_T=0$ (since it is locally free, being the branch locus of
$\pi$ contained in a normal crossings divisor (see \cite{Ko2})) and
$R^{2}g_{\ast }=0$ by reason of fibre dimension we have

\begin{equation}\label{eqdos}
R^{2}f_{\ast }({\cal M})=R^{1}g_{\ast }(R^{1}\pi _{\ast
}({\cal M}))=R^{1}g_{\ast }(R^{1}\pi _{\ast }{\cal O}_{T}\otimes{\cal L}).
\end{equation}

We observe that $R^2f_{\ast}({\cal M})$ is locally free. Indeed
${\cal M}=\pi^{\ast}{\cal L}$; since ${\cal L}$
is torsion and ${\cal L}^{\otimes i}_{|C}\not= {\cal O}_{C}$
for $1\leq i\leq n-1$, we
can consider the induced \'etale base change of $\pi $:
$$
\xymatrix{
\widehat{T}\ar[r]\ar[d]_{\widehat{\pi }}&T\ar[d]^{\pi }\\
\widehat{S}\ar[r]\ar[dr]_{\widehat{g}}&S\ar[d]^{g}\\
{}&B
}
$$
and get
\begin{equation}\label{eqtres}
R^j\widehat{f}_{\ast
}\omega_{\widehat{T}/B}=\displaystyle\mathop{\oplus}^{n-1}_{i=0} R^jf_{\ast
}(\omega_{T/B}\otimes{\cal M}^{\otimes i})
\end{equation}
$$R^j\widehat{f}_{\ast
}{\cal O}_{\widehat{T}}=\displaystyle\mathop{\oplus}^{n-1}_{i=0} R^jf_{\ast
}({\cal M}^{\otimes i}).$$

Hence $R^2f_{\ast}({\cal M})$ is locally free, being a subsheaf
of $R^2{\widehat f}_{\ast}{\cal O}_{{\widehat T}}$ (which is locally
free by relative duality and \cite{Ko2}).

Finally, remember that for fibrations of surfaces over curves we have
(cf. \cite{X1}; for this we do not need the fibration to be relatively
minimal):
$$\mbox{deg }(\pi _{t})_{\ast }\omega_{F_{t}/C_{t}}=
\chi ({\cal O}_{F})-\chi ({\cal O}_{C})\chi
({\cal O}_{G});\quad\mbox{deg }g_{\ast }\omega_{S/B}
=\chi ({\cal O}_{S})-\chi ({\cal O}_{C})\chi ({\cal O}_{B}).$$

Now we can compute

$$\begin{array}{lll}
\chi ({\cal O}_{T})&=\chi(\pi _{\ast }{\cal O}_{T})-\chi(R^{1}\pi_{\ast }
{\cal O}_{T}) &\\
&=\chi ({\cal O}_{S})-\chi(R^{1}\pi _{\ast }{\cal O}_{T}\otimes{\cal L})
&\mbox{since}\quad {\cal L}\in  {\rm Pic}^0(S) \\
&=\chi ({\cal O}_{S})+\chi(R^{1}g_{\ast }(R^{1}\pi _{\ast
}{\cal O}_{T}\otimes{\cal L}))&\mbox{by Leray}\\
&=\chi ({\cal O}_{S})+\chi(R^{2}f_{\ast }({\cal M}))&\mbox{by (\ref{eqdos})}\\
&=\chi ({\cal O}_{S})+\chi(f_{\ast }(\omega_{T/B}\otimes{\cal M}^{-1})^{\ast})&\mbox{by
relative duality}\\
&=\chi ({\cal O}_{S})-\mbox{deg}(f_{\ast
}(\omega_{T/B}\otimes{\cal M}^{-1}))+h^{0}(F,\omega_F \otimes {\cal M}_
{\vert F}^{-1})\chi({\cal O}_B)&
\mbox{by R.R. on $B$}
\end{array}
$$

By the choice of $\cal L$, Serre duality on $B$ and relative
duality on $B$ and $C$, we obtain

\begin{equation}\label{eqcinc}
h^0(F,\omega_{F}\otimes{\cal M}^{-1}_{|F})=
h^0(C,\pi _{\ast }\omega_{F/C}\otimes \omega_{C}\otimes {\cal L}^{-1}_{|C})=
\chi(\pi _{\ast }\omega_{F/C}\otimes \omega_{C}\otimes{\cal L}^{-1}_{|C})
\end{equation}
\begin{eqnarray*}
&&=-\chi(R^{1}\pi _{\ast }{\cal O}_{F}\otimes{\cal L})=-\chi(R^{1}\pi _{\ast
}{\cal O}_{F})=(\chi ({\cal O}_{F})-
\chi ({\cal O}_{C})\chi ({\cal O}_{G}))-g(G)\chi ({\cal O}_{C})\\
&&=\chi ({\cal O}_{F})-\chi ({\cal O}_{C})
\end{eqnarray*}

Hence:

\begin{equation}\label{eqsis}
\chi_f=\chi ({\cal O}_{F})\chi ({\cal O}_{B})-\chi
({\cal O}_{T})=\mbox{deg}(f_{\ast }(\omega_{T/B}\otimes{\cal M}^{-1}))-\mbox{deg}
(g_{\ast}\omega_{S/B})
\end{equation}

Now we use the hypothesis: $\pi _{\ast }\omega_{T/S}$ has a section and
hence we have an injection
$$
0\longrightarrow {\cal O}_{S}\longrightarrow \pi _{\ast }\omega_{T/S}
$$
which gives
$$
0\longrightarrow \omega_{S/B}\otimes{\cal L}^{-1}\longrightarrow \pi _{\ast
}\omega_{T/S}\otimes \omega_{S/B}\otimes{\cal L}^{-1}=\pi _{\ast
}(\omega_{T/B}\otimes{\cal M}^{-1})
$$
and so

\begin{equation}\label{eqset}
0\longrightarrow g_{\ast
}(\omega_{S/B}\otimes{\cal L}^{-1})\buildrel \tau \over
\longrightarrow g_{\ast }(\pi _{\ast
}(\omega_{T/B}\otimes{\cal M}^{-1}))=f_{\ast }(\omega_{T/B}\otimes{\cal M}^{-1})
\end{equation}

Note that deg$(g_{\ast }(\omega_{S/B}\otimes{\cal L}^{-1}))=
\mbox{deg}g_{\ast}\omega_{S/B}$. Indeed, by the choice of ${\cal L}$ we
have that $h^1(C,\omega_C\otimes {\cal L}_{\vert C}^{-1})=0$; on the
other hand $R^{1}g_{\ast }(\omega_{S/B}\otimes{\cal L}^{-1})$ is
locally free, being a subsheaf of $R^1{\widetilde g}_{\ast}
\omega_{{\widetilde S}/B}$ for the \'etale cover ${\widetilde S}
\longrightarrow S$ induced by ${\cal L}$. Hence
$R^{1}g_{\ast }(\omega_{S/B}\otimes{\cal L}^{-1})=0$.
Since ${\rm deg}R^{1}g_{\ast}(\omega_{S/B})={\rm deg}({\cal O}_S)=0$
we obtain the desired result using that
$\sum\limits_{j=0}^1(-1)^j\mbox{deg}R^jg_{\ast}(\omega_{S/B} \otimes
{\cal L})$ is independent of ${\cal L} \in {\rm Pic}^{0}(S)$.

In order to finish the proof, it suffices to check that $f_{\ast
}(\omega_{T/B}\otimes{\cal M}^{-1})$ is nef.
By (2) we have that
$f_{\ast
}(\omega_{T/B}\otimes{\cal M}^{-1})=f_{\ast }(\omega_{T/B}\otimes {\cal M}^{\otimes(n-1)})$ is
nef since it is a quotient of a nef vector bundle.

\def\circB{{\displaystyle\mathop{B}^{\circ}}}
\def\circT{{\displaystyle\mathop{T}^{\circ}}}
\def\circS{{\displaystyle\mathop{S}^{\circ}}}

\smallskip

(iii) Assume that for general  $t\in B$ we have a fibration
$h_{t}:F_{t}\longrightarrow C_{t}$. Let
$\circB\subseteq B$ be a non-empty open set such that
$f^{\circ}:\circT\longrightarrow \circB$ is smooth and for every
$t\in\circB$ there exists such a $h_{t}$. We can now consider the
fibration of abelian varieties $\psi
:\mbox{Alb}_{\circT/\circB}\longrightarrow \circB$. For every
$t\in\circB$ we have an abelian subvariety
$K_{t}=\mbox{ker}(\mbox{Alb}F_{t}\longrightarrow
\mbox{Alb}C_{t})\hookrightarrow \mbox{Alb}F_{t}=\psi ^{-1}(t)$. Then we
can apply \cite{BN} Theorem 2.5 and get, after a base change, a relative
abelian subvariety $K\hookrightarrow \mbox{Alb}_{\circT/\circB}$ over
$\circB$. Let $J=\mbox{Alb}_{\circT/\circB}\Big /K$. Consider the
natural map, after a base change,  $\varphi :\circT\longrightarrow
\mbox{Alb}_{{\circT/\circB}}\longrightarrow J$ over $\circB$. For
general $t\in\circB$, $\varphi _{t}:F_{t}\longrightarrow J_{t}$ has as
its image $C_{t}$  by construction. Let $\circS=\varphi (\circT)$ and
complete the map to get
$$
\xymatrix{
{}&\widetilde{T}\ar[d]\ar[dr]^{\pi }&{}\\
{T}\ar[d]_{f}&\overline{T}\ar[l]_{\bar\sigma}
\ar[r]\ar[d]_{\bar f}&{S}\ar[dl]\\
B&\overline{B}\ar[l]^{\sigma }&
}
$$
Note that we are in the same situation for $\bar f$ as in (ii). We have
even more since by construction the hypothesis $h^0(\pi_{\ast}\omega_
{T/S}) >0$ holds; indeed, let $E$ be a general curve on $S$ and let
$H$ be its pullback on $T$. We have that

$$\pi_{\ast}\omega_{H/E}=(\pi_{\ast}\omega_{T/S})\otimes {\cal O}_{E}.$$

If $E$ is ample enough, we also have

$$h^{0}(S,\pi_{\ast}\omega_{T/S}\otimes{\cal O}_{S}(-E))=
h^{1}(S,\pi_{\ast}\omega_{T/S}\otimes{\cal O}_{S}(-E))=0.$$

Hence $h^{0}(S,\pi_{\ast}\omega_{T/S})=h^{0}(E,\pi_{\ast}\omega_{H/E})$.
But in the case of fibred surfaces, $h^{0}(E,\pi_{\ast}\omega_{H/E}) \geq
q(H)-g(E)$ holds according to Fujita's decomposition (see Proposition 1.8
(i)). Finally
note that by \cite{Cai} $q(H)-g(E) \geq q(T)-q(S) \geq 1$.

So we can apply the same argument as in (ii) and get
$\chi_{\bar f}\geq 0$.

Then we have
$$
\xymatrix{
\overline{T}\ar[r]^{\alpha }\ar[dr]_{\bar
f}&T\displaystyle\mathop{\times
}_{B}\overline{B}\ar[d]^{f'}\ar[r]&T\ar[d]^{f}\\
{}&\overline{B}\ar[r]_{\sigma }&B }
$$
where  $\alpha $ is induced by $\overline{f}$ and $\overline{\sigma}$. Since
$f$ is semistable we can apply base change theorem (\cite{M}, 4.10)
and get
$$
\overline{f}_{\ast }\omega_{\overline{T}/\overline{B}}=f'_{\ast
}\omega_{{T{\displaystyle\mathop{\times }_{B}}\overline{B}}/\overline{B}}=\sigma
^{\ast }(f_{\ast }\omega_{T/B})\, .
$$
In fact we also have the same equality for $R^{1}f_{\ast }$: take
${\cal H}$ a very ample line bundle in $T$ and let $H$ be a general smooth
member of its associated linear system. We have in a natural way
$$
0\longrightarrow f_{\ast }\omega_{T/B}\longrightarrow f_{\ast
}(\omega_{T/B}\otimes {\cal H})\longrightarrow f_{\ast }\omega_{H/B}\longrightarrow
R^{1}f_{\ast }\omega_{T/B}\longrightarrow 0
$$
since  $R^{1}f_{\ast }(\omega_{T/B}\otimes{\cal H})=0$
(by Kodaira vanishing $h^1(F,\omega_F \otimes {\cal H}_{\vert F})=0$
and $R^1f_{\ast}(\omega_{T/B}\otimes {\cal H})$ is locally free by the
trick of a cyclic cover used in (ii)).

Note that all of them are locally free. Hence we have that after taking
$\sigma^{\ast}$ we still have a long exact sequence. Considering
the analogous exact sequence for $\overline f$ and the natural maps we get
$$
\xymatrix{
0\ar[r]&\sigma ^{\ast }(f_{\ast }\omega_{T/B})\ar[r]&\sigma ^{\ast
}(f_{\ast }\omega_{T/B}\otimes{\cal H})\ar[r]&\sigma ^{\ast }(f_{\ast }\omega_{H/B})
\ar[r]&\sigma ^{\ast }(R^{1}f_{\ast }\omega_{T/B})\ar[r]&0\\
0\ar[r]&\bar f_{\ast }(\omega_{\overline T/\overline
B})\ar[r]\ar[u]^{\cong}&\bar
f_{\ast }(\omega_{\overline T/\overline B}\otimes{\overline {{\cal H}}})\ar[r]\ar[u]&\bar
f_{\ast } \omega_{\overline H/\overline B} \ar[r]\ar[u]^{\cong}&R^{1}\bar
f_{\ast }\omega_{\overline T/\overline B}\ar[r]\ar[u]^{\gamma }&0}
$$
where  $\gamma $ is naturally induced and exhaustive. Since
$R^{1}\bar f_{\ast }\omega_{\overline T/\overline B}$ and $\sigma ^{\ast
}(R^{1}f_{\ast }\omega_{T/B})$ are both locally free sheaves of the same rank
over $\overline B$, $\gamma $ is an isomorphism.

So we have
$$
0\leq \chi_{\bar f}=\mbox{deg}\bar f_{\ast }\omega_{\overline T/\overline
B }-\mbox{deg }R^{1}\bar f_{\ast }\omega_{\overline T/\overline
B}=n(\mbox{deg}f_{\ast }\omega_{T/B}-\mbox{deg}R^{1}f_{\ast }\omega_{T/B})
=n \chi_f \, .
$$

\smallskip

(iv) Since  $T$ is of Albanese general type, then so is  $F_{t}$ for
$t\in B$ general. We can apply then \cite{Be1} Theorem 1 to get that
$\{{\cal L} \in  {\rm Pic}^0(F_{t})|\, h^{1}(F_{t},{\cal L})\not= 0\}$ is the union of
subtori $h^{\ast }_{i}({\rm Pic}^0(C_{i}))$ for fibrations
$h_{i}:F_{t}\longrightarrow C_{i}$ with  $g(C_{i})\geq 2$ and a finite
number of (torsion) points.

Under our assumptions we can take an $n$-torsion
element ${\cal L}\in  {\rm Pic}^0(T)$ such that for $1\leq i\leq n-1$,\ $h^{1}
(F_{t},{\cal L}^{\otimes i}_{|F_{t}})=0 $.

Hence as in Lemma 1.5 (iii) if we consider the \'etale cover
$\sigma :\widetilde{T}\longrightarrow T$ associated to ${\cal L}$, and
$\widetilde{f}=f\circ\sigma $\ we have that $f_{\ast }(\omega_{T/B}\otimes{\cal L})$
is a quotient of $\widetilde{f}_{\ast }\omega_{\widetilde{T}/B}$ hence it is
nef.

Since $h^{1}(F_{t},{\cal L}_{|F_{t}}^{-1})=0$ we have  $R^{1}f_{\ast
}(\omega_{T/B}\otimes {\cal L})=0$ (as above it is locally free) and hence
\begin{eqnarray*}
\chi_f&=&\mbox{deg}f_{\ast }\omega_{T/B}-\mbox{deg}R^{1}f_{\ast
}\omega_{T/B}=\\
&=&\mbox{deg}f_{\ast }(\omega_{T/B}\otimes{\cal L})-\mbox{deg}R^{1} f_{\ast
}(\omega_{T/B}\otimes {\cal L})=\mbox{deg}f_{\ast }(\omega_{T/B}\otimes {\cal L})\geq 0
\end{eqnarray*}
\quad $\Box$

We finish with the following result on the structure of ${\cal E}
=f_{\ast }w_{T/B}$. The first part is a well known result of Fujita
(\cite{F1}, \cite{F2}). The second part can be found in \cite{BT}.
We include a brief idea of proof for benefit of the reader.

{\bf Proposition 1.8.} \quad
{\it Let $X$, $Y$ be smooth varieties of
dimensions  $n>m$ respectively. Let $f:X\longrightarrow Y$ be a fibration
with a simple normal crossings branching.

Then:
\begin{enumerate}
\item[{\rm (i)}] If $m=1$, ${\cal E} =f_{\ast
}w_{X/Y}=\A\oplus\displaystyle\mathop{\oplus}^{r}_{i=1}E_{i}\oplus
{\cal O}^{h}_{Y}$ where $\A$ is an ample vector bundle
(or zero),
$E_{i}$ are stable, degree zero non-trivial vector bundles and
$h=h^{1}(Y,f_{\ast}\omega_{X})$. If $X$ is a surface ($n=2$), then
$h=q(X)-q(Y)$.
If $E$ is a stable degree zero vector bundle which is a quotient of
${\cal E}$, then $E$ is one of the $E_i$ or ${\cal O}_{Y}$.

\item[{\rm (ii)}] If there exists a vector bundle $E$ with
${\rm det}(E)=\L\in{\rm Pic}^{\!\circ}(Y)$ and an epimorphism
${\cal E} =f_{\ast }w_{X/Y}\longrightarrow E$, then  $\L$ is
torsion. In particular, when  $m=1$, $E_{i}$ is a torsion line bundle
whenever rank $(E_{i})=1$.
\end{enumerate}
}

\dem
(i) When Y is a curve we have a decomposition $f_{\ast}\omega _{X/Y}=
{\cal A} \oplus {\cal U}$, where ${\cal A}$ is ample (or zero) and
${\cal U}$ is flat (see \cite {F2}). Since $Y$ is a curve and ${\cal U}$
is flat, we have a decomposition of ${\cal U}$ in direct sum of
stable, degree zero pieces (see, for example, \cite{Be1}). Finally
we can use \cite{F1}.

(ii) By induction on the dimension of $Y$ we can assume that $Y$ is a
curve. It its easy to see that $\L$ is a quotient of
$f^{(s)}_{\ast}\omega _{X^{(s)}/Y}$ (the s-th fibred product of $f$
over $Y$ (see \cite{V})). Then we can check that $\L$ increases the
general value of $h^{0}(Y,(R^{d}f_{\ast}^{(s)}{\cal O}^{(s)})\otimes
{\cal M})$, for ${\cal M} \in {\rm Pic}^{0}Y$. Then we can apply a
result of Simpson (see \cite{Si}) to get that $\L$ is torsion.
\quad $\Box$

\vglue.5truecm

\noindent
{\bf  2. Lower bounds for the slopes of fibred threefolds}

\medskip

\indent
We give here a lower bound for $\lambda_1(f)$ (and hence for $\lambda_2
(f)$ provided it is well defined) in the case of a relatively minimal
fibred threefold with fibres of general type. The bounds we obtain
are considerably better than Ohno's bounds (\cite{O} Main Theorem 1)
as long as $p_g(F)\gg0$.

First we need some results on linear systems on surfaces of general type.

{\bf Lemma 2.1.} \quad
{\it
Let $F$ be a minimal surface of general type
such that $p_{g}(F)\geq 3$ and let $\tau :\widetilde{F}\longrightarrow
F$ be a birational morphism.
Let $0 \leq P \leq Q \leq \tau^{\ast}K_F$ be two nef and effective
divisors, such that the complete linear systems $\vert P \vert$ and
$\vert Q \vert$ are base point free. Let $r \leq s$ be the dimensions of
$H^{0}(F,{\cal O}_F(P))$ and $H^{0}(F,{\cal O}_F(Q))$ respectively.
Let $\Sigma$ be
the image of $\widetilde F$ through the map $\varphi$ induced by
$\vert P \vert$. Then:
\begin{enumerate}
\item[{\rm (i)}]If $\varphi $ is a generically finite map, then we have
\begin{itemize}
\item  $P(\tau ^{\ast }K_{F})\geq P^{2}\geq 2r-4+2q(\Sigma )$ if
$\varphi $ is a double cover of a geometrically ruled surface $\Sigma$.
\item $P(\tau ^{\ast }K_{F})\geq P^{2}\geq 3r-7$  otherwise.
\end{itemize}
\item[{\rm (ii)}] If $\vert P\vert$ is composed with a pencil of curves $D$ of
(geometric) genus $g$,  $\widehat D=\tau _{\ast }D$, and  $\vert Q\vert$
induces a
generically finite map, then we have
\begin{itemize}
\item  $QP\geq 2(r-1)$.
\item  $QP\geq 3(r-1)$ except if $D$ is hyperelliptic and
$\vert Q\vert_{|D}=g^{1}_{2}$.
\item  $QP\geq 4(r-1)$ except if $D$ is hyperelliptic or trigonal and
$\vert Q\vert_{|D}=g^{1}_{2}$ or $g^{1}_{3}$.
\item  $QP\geq 5(r-1)$ except if $D$ is hyperelliptic, trigonal or
tetragonal and $\vert Q \vert_{\vert D}=g_2^1$, $2g_2^1$, $g_3^1$ or
$g_4^1$, or $D$ is of genus $2$ or $3$.
\item  $P(\tau^{\ast }K_{F})\geq (2g-2)(r-1)$ or $(2g-2)r$, according to
whether the pencil is rational or not, if the pencil  $|\widehat D|$ in $F$ has
no base point.
\item $P(\tau ^{\ast }K_{F})\geq (2p_{a}(\widehat D)-2-\widehat
D^{2})(r-1)$ if the pencil $|\widehat D|$ has some base point.
\end{itemize}
\item[{\rm (iii)}] If $|K_{F}|$ is composed with a pencil which general member
$D$ is as
in (ii), then we have
\begin{itemize}
\item $P(\tau ^{\ast }K_{F})\geq (2g-2)(r-1)$ or $(2g-2)r$, according to
whether the pencil is rational or not, if the pencil $|\widehat D|$ in $F$ has
no base point.

\item  $P(\tau ^{\ast }K_{F})\geq \mbox{\rm max}\{ \sqrt{2(g-1)\left
(1-\frac{1}{p_{g}(F)} \right )(p_{g}(F)-1)\widehat D^{2}}$, \break
$(2p_a(\widehat D)-2-\widehat D^2)(r-1)\}$ otherwise.
\end{itemize}
\end{enumerate}
}

\dem
(i) Since  $P$ is nef and $P \leq \tau ^{\ast }K_{F}$, we obviously
have $P(\tau ^{\ast }K_{F})\geq P^{2}$. It is a well known fact that deg
$\Sigma \geq r-2+q(\Sigma )$ if $\Sigma $  is geometrically ruled and
that deg $\Sigma \geq 2r-4$ otherwise (see \cite{Be2}).

Let $a=\mbox{deg }\varphi $. If $a\geq 3$ then $P^{2}\geq
3\mbox{deg}\Sigma \geq 3(r-2)>3r-7$. If  $a=2$ and $\Sigma $ is not
geometrically ruled, then $P^{2}\geq 2(2r-4)=4r-8\geq 3r-7$. If $\Sigma $ is
geometrically ruled, then $P^{2}\geq 2\mbox{deg}\Sigma \geq 2r-4+2q(\Sigma )$.

If $a=1$, let $C\in \vert P\vert $ be a smooth curve ($\vert P\vert $ has
no base point). Then
$2P_{|C}\leq (\tau ^{\ast }K_{F}+P)_{|C}\leq
(K_{\widetilde{F}}+P)_{|C}=K_{C}$. So  $\mbox{deg}P_{|C}\leq g(C)-1$. We
can then apply \lq\lq Clifford plus\rq\rq \ lemma (cf. \cite{Be4})
and get $P^{2}=\mbox{deg}P_{|C}\geq 3h^{0}(C,P_{|C})-4\geq
3h^{0}(\widetilde{F},P)-7=3r-7$.

\vglue.5truecm

(ii) Let $\varphi_P(F)=C \subseteq {\proj}^{r-1}$. The map $F
\longrightarrow C$ may not have connected fibres; consider the Stein
factorization of $\varphi_P$, $F \longrightarrow \widetilde C
\longrightarrow C$. Note that then we have $P \equiv \alpha D$ where
$D$ is an irreducible smooth curve such that $D^2 =0$ and
$\alpha=\alpha_1 \alpha_2$ where $\alpha_1=\mbox{deg}(\widetilde C
\longrightarrow C)$ and $\alpha_2 \geq r-1$ (and equality holds only
when $C$ is rational). The pencil $\vert P \vert$ is said to be rational
if $\widetilde C={\proj}^1$ and irrational otherwise. Note that in
general also $\alpha \geq r-1$ and $\alpha \geq r$ if the pencil is
irrational.

Since $P \leq Q$, the map $\varphi_P$ factors through $\varphi_Q$.
Let $\Sigma=\varphi_Q(F)$ and consider the induced map $\psi:\Sigma
\longrightarrow C$.
By construction clearly $\varphi_Q(D)\subseteq \psi^{-1}(t)$ for
$t \in C$ (note that $\psi^{-1}(t)$ does not need to be connected).

We have that $QP=\alpha_2(\alpha_1 QD)\geq (r-1)(\alpha_1QD)$. Let $a$
be the degree of $\varphi_{Q \vert D}$, $\overline D=\varphi_Q(D)$ and
$d=\mbox{deg}\overline D$.  Note that $a$ divides deg$\varphi_Q$
although we will not use it. We have then $\alpha_1QD=\alpha_1 ad$.
Note that $ad \geq 2$ (otherwise $\widetilde F$ would be covered by rational
curves) and hence $QP \geq 2(r-1)$. But if $QP<3(r-1)$, then $\alpha_1=1$,
$a=2$, $d=1$ (if $a=1$, $d=2$ again $\widetilde F$ is covered by
rational curves). Hence $D$ is hyperelliptic and $\vert Q\vert_{\vert D}
=g^1_2$.

If $QP<4(r-1)$, then $\alpha_1ad \leq 3$. If $\alpha_1ad=2$, the previous
argument holds. If $\alpha_1ad = 3$, then $\alpha_1=1$, $a=3$, $d=1$
(if $a=1$, $d=3$, $\widetilde F$ would be covered by elliptic or rational curves,
a contradiction since $F$ is of general type). Then $D$ is trigonal and
$\vert Q \vert_{\vert D}=g_3^1$.

If $QP<5(r-1)$, then $\alpha_1ad \leq 4$ and we must only study the case
$\alpha_1ad=4$. Four possibilities may occur. Either $\alpha_1=2$,
$a=2$, $d=1$ (then $D$ is hyperelliptic and $\vert Q \vert _{\vert D}
=g_2^1$) or $\alpha_1=1$, $a=4$, $d=1$ (then $D$ is tetragonal and
$\vert Q \vert_D =g_4^1$) or $\alpha_1=1$, $a=2$, $d=2$ (then $D$
is again hyperelliptic and $\vert Q \vert_D =2g_2^1$) or
$\alpha_1=a=1, d=4$ (then $D$ has at most genus 3; in particular
$D$ is also hyperelliptic or trigonal).

Finally note that $D\tau^{\ast}K_F=2g-2$ if $\vert \widehat D \vert$ has
no base point and $D \tau^{\ast}K_F=\widehat DK_F=2p_a(\widehat D)
-2-\widehat D^2$ (by adjunction formula on $F$) otherwise. Hence the result
follows from $P(\tau^{\ast}K_F)=\alpha D(\tau^{\ast}K_F)$ and the
previous bound of $\alpha$.

\vglue.5truecm

(iii) The first result is analogous to (ii) keeping in mind that if $P\leq
\tau ^{\ast }K_{F}$, then  $\vert P\vert$ is composed with a pencil of the same
genus as $|K_{F}|$ and with the same general fibre.

Part of the second statement follows as in (ii). For the rest recall that
from Hodge Index theorem  $(K_{F}\widehat
D)^{2}\geq K^{2}_{F}\widehat D^{2}$ and that when $|K_{F}|$  is composed
with a pencil of genus zero, then
$K^{2}_{F}\geq 2(g-1)$
$\left
(1-\frac{1}{p_{g}(F)} \right )
(p_{g}(F)-1)$ (\cite{K6}, Lemma 3.3).
\quad $\Box$

{\bf Definition 2.2.} \quad
We say that a linear pencil $\vert Q\vert$ on $F$
is of type $(r,g,p)$, $r\geq 2$, $g\geq 2$, $p=0,1$ if
$\vert Q\vert$ is a {\rm complete} linear system of  $r$-gonal
(but not $s$-gonal for $s<r$) curves of (geometric) genus $g$,
rational if $p=0$, irrational if $p=1$.
If $\vert Q \vert$ is a rational pencil, we call ${\widehat D}$
a generic member and $\delta=K_F{\widehat D}=2p_a({\widehat D})
-2-{\widehat D}^2$. If $\vert Q \vert$ is base point free, then
clearly $\delta=2g-2$.

{\bf Remark 2.3.} \quad
We recall the following result, originally
due to  Xiao (\cite{X1}) for the case of surfaces and to  Ohno
(\cite{O}) for 3-folds and given in full generality by Konno in
\cite{K1}. See this last paper for details.

Given  $f:T\longrightarrow B$ as before, consider  ${\cal E}
=f_{\ast }w_{T/B}$ and the Harder-Narasimhan filtration of
${\cal E}$:
$$
0={\cal E}_0 \subseteq \dots \subseteq {\cal E}_{\ell -1} \subseteq
{\cal E}_{\ell}={\cal E}
$$
and $\mu _{1}\!>\!\mu _{2}\!>\!\dots\!>\!\mu _{\ell}\!\geq \!0$, where $\mu _{i}=\mu
({\cal E} _{i}/{\cal E} _{i-1})$. If  $r_{i}={\rm rk}{\cal E} _{i}$, then
$\mbox{deg }{\cal E}
=\displaystyle\sum^{\ell}_{i=1}r_{i}(\mu _{i}-\mu _{i-1})$.

Consider the induced rational maps $\varphi _{i}:T ---> {\proj}_
{B}({\cal E} _{i})$ and let
$\widetilde{T}\displaystyle\mathop{\longrightarrow }_{\alpha }T$ be a
desingularization of $T$ which eliminates the indeterminacy of the
$\varphi _{i}$. Let  $\tau =\alpha
_{|\widetilde{F}}:\widetilde{F}\longrightarrow F$ the restriction of
$\alpha $ to a general fibre $\widetilde{F}$ of $f\circ\alpha$. Then we have
a sequence of nef effective Cartier divisors on ${\widetilde F}$,
$P_{1}<P_{2}<\dots<P_{\ell}\leq P_{\ell+1}\leq \tau ^{\ast }K_{F}$
(where $P_{\ell +1}$ can be chosen to be $P_{\ell}$ or $\tau^{\ast}
K_F$) such that the (projective) dimension of the linear systems
$\vert P_i \vert$ are $r_i-1$
($r_{i}=\mbox{rk}({\cal E} _{i})$) and such that
$(P_{i}-\mu _{i}{\widetilde F})$ is nef and
for any $1\leq
i<\dots<i_{n}\leq i_{n+1}=\ell+1$ and any $1\leq m\leq n$, we get
\begin{eqnarray}
K^{3}_{T/B}&=&(\mu ^{\ast }K_{T/B})^{3}\geq \sum^{m-1}_{p=1}
(P_{i_{p}}+P_{i_{p+1}}) P_{i_{m}} (\mu _{i_{p}}-\mu _{i_{p+1}}) +
\label{eqdosuna}
\\
&&\quad +
\sum^{n}_{p=m} (P^{2}_{i_{p}}+P_{i_{p}}P_{i_{p+1}}+ P^{2}_{i_{p+1}})
(\mu _{i_{p}}-\mu _{i_{p+1}})\nonumber\\
K^{3}_{T/B}&\geq &\sum^{n}_{p=1} (P_{i_{p}}+P_{i_{p+1}}) (\tau
^{\ast }K_{F}) (\mu _{i_{p}}-\mu _{i_{p+1}})\label{eqdosdos}
\end{eqnarray}
where $\mu _{\ell +1}=0$.

In particular, for the indices $\{1,\ell\}$, we have
\begin{equation} \label{eqdostres}
K^{3}_{T/B}\geq P^{2}_{\ell}\mu _{1}+(P^{2}_{\ell}+P_{\ell +1}^2)\mu _{\ell}
\geq P^{2}_{\ell} (\mu _{1}+2\mu _{\ell}).
\end{equation}

Now we can state the main result of this section:

{\bf Theorem 2.4.} \quad
{\it
Let $T$ be a normal,
projective threefold with only canonical singularities and let
$f:T\longrightarrow B$ be a relatively minimal fibration onto a smooth
curve of genus $b$. Let $F$ be a general fibre.

Assume $F$ is of general type, $p_{g}(F)\geq 3$ and that $\chi
_{f}=\chi ({\cal O}_{F})\chi ({\cal O}_{B})-\chi ({\cal O}_{T})>0$. Then:
\begin{enumerate}
\item[{\rm (i)}] If $|K_{F}|$ is not composed with a pencil and $|K_{F}|$ has
no subpencil $\vert P\vert$ of type  $(r,g,p)$, $r=2,3,4$, then
$$
\lambda _{2}(f)\geq \lambda _{1}(f)\geq 9\left  (1-\frac{17}{3p_{g}(F)
+10} \right )$$
\item[{\rm (ii)}] If $|K_{F}|$ is composed with a pencil with generic fibre
$\widehat D$ of (geometric) genus $g$, then
$$
\begin{array}{ll}
\lambda _{2}(f)\geq \lambda _{1}(f)\geq 4g-4&\mbox{if the pencil is irrational}\\
\lambda _{2}(f)\geq \lambda _{1}(f)\geq
\frac{p_{g}(F)}{p_{g}(F)+1}(4p_{a}(\widehat D)-4-2\widehat
D^{2})&\mbox{if the pencil is rational}
\end{array}
$$
\item[{\rm (iii)}] If $|K_{F}|$ is not composed, has no hyperelliptic
subpencils, and has  subpencils
$\{\vert Q_{i}\vert \}_{i \in I}$ of type $(r_{i},g_{i},p_{i})$ with
$r_{i}=3$ or $4$ then
$$
\lambda _{2}(f)\geq \lambda _{1}(f) \geq \displaystyle \min_{i \in I} \{\lambda
_{r_{i}}^{p_{i}}(Q_{i})\},
$$
where

$\qquad
\lambda^{1}_{3}(Q)=9-\frac{9}{4g-7}-{\varepsilon}^{1}_{3}(p_{g}(F),g)$

\smallskip

$\qquad
\lambda^{0}_{3}(Q)=\left \{
\begin{array}{ll}
9-\frac{9}{\delta -3}-{\varepsilon}^{0}_{3}(p_{g}(F),\delta)
&\quad \mbox{if}\quad \delta \geq 7,\\
6\left (1-\frac{17}{3p_g(F)+10}\right ) &\quad \mbox{otherwise}
\end{array} \right .$

\smallskip

$\qquad
\lambda^{1}_{4}(Q)=9-\frac{3}{4g-9}-{\varepsilon}^{1}_{4}(p_{g}(F),g)$

\smallskip

$\qquad
\lambda^{0}_{4}(Q)=\left \{\begin{array}{ll}
9-\frac{3}{\delta -5}-{\varepsilon}^{0}_{4}(p_{g}(F),\delta)
&\quad \mbox{if}\quad \delta \geq 7, \\
8\left (1-\frac{17}{3p_g(F)+10}\right )&\quad \mbox{otherwise}
\end{array}\right .$

\item[{\rm (iv)}] If $\vert K_{F} \vert$ has hyperelliptic subpencils
$\{\vert Q_{i}\vert \}_{i \in I}$ of type $(2,g_{i},p_{i})$ then

$$
\lambda _{2}(f)\geq \lambda _{1}(f) \geq \displaystyle \min_{i \in I} \{\lambda
_{r_{i}}^{p_{i}}(Q_{i})\}
$$

where

$\qquad
\lambda^{1}_{2}(Q)=
\left \{\begin{array}{ll}
6-\frac{2}{2g-3}-{\varepsilon}^{1}_{2}(p_{g}(F),g)
&\quad \mbox{if}\quad g \geq 3, \\
4\left (1-\frac{1}{p_g(F)-1}\right )
&\quad \mbox{if} \quad g=2 \end{array} \right .$

\smallskip

$\qquad
\lambda^{0}_{2}(Q)=
\left \{\begin{array}{ll}
6-\frac{4}{\delta -2}-{\varepsilon}^{0}_{2}(p_{g}(F),\delta)
&\quad \mbox{if}\quad \delta \geq 4, \\
4 \left (1-\frac{9}{2p_g(F)+5}\right )
&\quad \mbox{otherwise} \end{array}\right .$

\end{enumerate}
and where ${\varepsilon}^{p}_{r} \sim O(\frac{1}{p_g(F)})$ are the following
positive functions

{\scriptsize
$$
\varepsilon^1_4=\frac{(68g-159)(36g-84)}{(4g-9)^2(3p_g(F)-7)+
(68g-159)(4g-9)} \quad
\varepsilon^0_4=\frac{(17\delta -91)(9\delta -48)}{(\delta -5)^2
(3p_g(F)-7)+(\delta -5)(17 \delta -91)}$$
$$\varepsilon^1_3=\frac{36(g-2)(68g-137)}{(4g-7)^2(3p_g(F)-7)+
(68g-137)(4g-7)} \quad
\varepsilon^0_3=\frac{(9 \delta -36)(17 \delta -69)}{(\delta-3)^2
(3p_g(F)-7)+(\delta -3)(17 \delta-69)}$$
$$\varepsilon^1_2=\frac{4(3g-5)(9g-17)}{(2g-3)^2(p_g(F)-2)+
(9g-17)(2g-3)} \quad \quad
\varepsilon^0_2=\frac{(9\delta -24)(6 \delta -16)}{(\delta -2)^2
(2p_g(F)-4)+(\delta -2)(9 \delta -24)}$$
}
}

{\bf Remark 2.5.} \quad
The statement of Theorem 2.4 looks considerably simpler if we
look at its asymtotic behaviour as $p_g(F)$ tends to infinity
(in which case the functions $\varepsilon$
tend to be zero). This behaviour will play a special role in the next
section. We also observe that even if $\chi_f <0$, then the bounds in
the theorem hold for $\lambda_1(f)$ as far as $\Delta_{f} \neq 0$.

\dem We consider ${\cal E}=f_{\ast} \omega_{T/B}$ and its
Harder-Narasimhan filtration
$$0 ={\cal E}_0\subseteq {\cal E}_1 \subseteq \dots \subseteq
{\cal E}_\ell ={\cal E}$$
with slopes $\mu_1 > \mu_2 > \dots > \mu_\ell \geq 0$ and ranks
$1 \leq r_1 <r_2 < \dots < r_\ell = p_g(F)$. As in \S 1.2, each piece
induces a Cartier divisor $P_i$ on $F$ such that the linear system
$\vert P_i \vert $ has projective dimension at least $r_i-1$. We define
as usual $\mu_{\ell +1}=0$, $P_{\ell +1}=P_\ell$.
We observe that following Remark 2.3 we could define, if necessary, $P_{\ell +1}=\tau^{\ast}
K_F \geq P_{\ell}$ (where $\vert P_{\ell} \vert$ is the moving part of
$\vert K_F \vert$), although this possibility will only be used
in very special computations and will be specifically pointed out.
Remember that
we have $\Delta_f=\mbox{deg}{\cal E}=\sum\limits_{i=1}^\ell r_i (\mu_i -
\mu_{i+1})$.

\smallskip

Consider first the case where $|K_{F}|$ is composed. Using Lemma 2.1
(iii),
Remark 2.3 and that
$r_{i+1} \geq r_i+1$ we get, if the pencil is irrational
\begin{eqnarray*}
K^{3}_{T/B}&\geq&\sum^{\ell}_{i=1}(P_{i}+P_{i+1})(\tau ^{\ast
}K_{F})(\mu _{i}-\mu _{i+1}) \\
&\geq& \sum^{\ell-1}_{i=1}((4g-4)r_{i}+(2g-2))(\mu _{i}-\mu _{i+1})+
(4g-4) r_{\ell}\mu _{\ell}\\
&=&(4g-4)\Delta_f+(2g-2)(\mu _{1}-\mu _{\ell})\geq (4g-4)\Delta_f.
\end{eqnarray*}

If the pencil is rational and $\widehat D$ is a generic member of its
moving part, then we have
\begin{eqnarray*}
K^{3}_{T/B}&\geq &\sum^{\ell-1}_{i=1}((4p_{a}(\widehat D)-4-2\widehat
D^{2})r_{i}-(2p_{a}(\widehat D)-2-\widehat D^{2}))(\mu _{i}-\mu
_{i+1})\\
&&\quad + (4p_{a}(\widehat D)-4-2\widehat D^{2})(r_{\ell}-1)\mu
_{\ell}\\
&=&(4p_{a}(\widehat D)-4-2\widehat D^{2})\Delta_f-(2p_{a}(\widehat
D)-2-\widehat D^{2})(\mu _{1}+\mu _{\ell}).
\end{eqnarray*}

By Remark 2.3 using the indices $\{1,\ell\}$, we get
\begin{eqnarray*}
K^{3}_{T/B}&\geq &(P_{1}+P_{\ell})(\tau ^{\ast }K_{F})(\mu _{1}-\mu
_{\ell})+2P_{\ell}(\tau ^{\ast }K_{F})\mu _{\ell} \\
&\geq &P_{\ell}(\tau ^{\ast }K_{F})(\mu _{1}+\mu _{\ell}) \\
&\geq & (2p_{a}(\widehat D)-2-\widehat D^{2})p_{g}(F) (\mu _{1}+\mu
_{\ell}).
\end{eqnarray*}

And hence eliminating $(\mu_1+\mu_\ell)$ from the above inequalities,
we get
$$
\left  (1+\frac{1}{p_{g}(F)} \right)K^{3}_{T/B}\geq
(4p_{a}(\widehat D)-4-2\widehat D^{2})\Delta_f
$$
which proves (ii).

\vglue.5truecm

From now on we assume that $|K_{F}|$ is not composed. Let
$$
m=\min \{ k\, |\,\, \vert P_{k} \vert \quad \mbox{induces a generically
finite map}\}
\leq \ell.
$$

By Remark 2.3 we have
$$
K^{3}_{T/B}\geq \sum^{m-1}_{i=1}(P_{i}+P_{i+1})P_{m}(\mu _{i}- \mu
_{i+1})+\sum^{n}_{i=m} (P^{2}_{i}+P_{i}P_{i+1} + P^{2}_{i+1}) (\mu
_{i}-\mu _{i+1}).
$$

Note that, for $i \geq m$, we have $P_{i+1}^2\geq P_iP_{i+1}\geq P_i^2$
and, if $P_iP_{i+1}=P_i^2$, then $P_i=P_{i+1}$. Indeed, we have $P_{i+1}
=P_i+D_i$, with $D_i \geq 0$. Hence $P_{i+1}^2=P_{i+1}(P_i+D_i)\geq
P_{i+1}P_i=(P_i+D_i)P_i\geq P_i^2$ since $P_i$ and $P_{i+1}$ are nef
and $D_i$ effective.
If $P_iP_{i+1}=P_i^2$, we would have $P_iD_i=0$. Since $\vert P_i\vert$ is
base point free and is not composed, Hodge Index Theorem applies and
hence either $D_i^2<0$ (which is impossible since then
$P_{i+1}^2=P_i^2+2P_iD_i+D_i^2<P_i^2$) or $D_i=0$. So we get $P_i=
P_{i+1}$.

Assume $\vert K_F \vert$ has no hyperelliptic subpencil
(in particular, the maps induced by the linear systems $\vert P_i \vert$
are never double covers of geometrically ruled surfaces). Then:
\begin{equation}\label{eqdosquatre}
\begin{array}{lll}
\mbox{for}&m\leq i\leq \ell-1,&P^{2}_{i}+P_{i}P_{i+1}+P^{2}_{i+1}\geq 9r_{i}-17\\
\mbox{and for}&i=\ell,&P^{2}_{\ell}+P_{\ell}P_{\ell+1}+P^{2}_{\ell+1}\geq
9r_{\ell}-21.
\end{array}
\end{equation}

Indeed, we denote by $\varphi_i$ the map induced by $\vert P_i \vert$ and
put $a_i=\mbox{deg}\varphi_i$. Note that $r_i \geq 3$. First consider
the case $m \leq i \leq \ell-1$.
By Lemma 2.1, if $\varphi_{i}$ and $\varphi_{i+1}$ are not double covers
of geometrically ruled surfaces, we have
$P^{2}_{i}\geq 3r_{i}-7$ and
$P^{2}_{i+1}\geq 3r_{i+1}-7\geq 3r_{i}-4$;
if $P_i \not= P_{i+1}$, then $P_{i}
P_{i+1}>P^{2}_{i}\geq 3r_{i}-7$ and we are done. If $P_i=P_{i+1}$, then
$(P_i^2+P_iP_{i+1}+P_{i+1}^2)=3P_{i+1}^2\geq 9 r_{i+1}-21 \geq
9r_i-12$.


If $i = \ell$ the result follows immediately from the previous
considerations.

Similarly, if $\vert K_F \vert$ admits hyperelliptic subpencils,
\begin{equation}\label{eqdoscinc}
\begin{array}{lll}
\mbox{for}&m\leq i\leq \ell-1,&P^{2}_{i}+P_{i}P_{i+1}+P^{2}_{i+1}\geq
6r_{i}-9\\
\mbox{and for}&i=\ell,&P^{2}_{\ell}+P_{\ell}P_{\ell+1}+P^{2}_{\ell+1}\geq
6r_{\ell}-12.
\end{array}
\end{equation}

Indeed, $9r_i-17 \geq 6r_i -9$ (since $r_i \geq 3$) so we only have to
check the case $a_i=2$ and $\varphi_i(\widetilde F)$ a geometrically
ruled surface. Since $\varphi_i$ factorizes through $\varphi_{i+1}$,
$a_{i+1}=1$ or 2. In any case $P_{i+1}^2\geq 2r_{i+1}-4\geq2r_i-2$. If
$P_i \not= P_{i+1}$, we have $P_i P_{i+1} > P_i^2 \geq 2r_i-4$
and we are done. If $P_i=P_{i+1}$,
then $P_i^2+P_{i+1}P_i+P_{i+1}^2=3P_{i+1}^2\geq
6r_i-6$. For $i=\ell$ the assertion is clear.

Observe that since $P_1 \leq \dots \leq P_{m-1}$, all the maps
$\varphi_i$ induced by $\vert P_i \vert$ are composed of the same
pencil (with the only exception of $r_1=1$, $P_1=0$ for which we
have no defined map $\varphi_1$). Indeed if
$i<j\leq m-1$ the map $\varphi_i$ factors through the map
$\varphi_j$. Since $\varphi_i(\widetilde F)$ and $\varphi_j(\widetilde
F)$ are curves, both maps have, after the Stein factorization, the same
fibre.

Let $P_{m-1}$ (and hence $P_i$ for $i \leq m-1$) be of type
$(r,g,p)$. Now, if $i=m-1$, then:
\begin{equation}\label{eqdossis}
(P_{m-1}+P_{m})P_{m}\geq
\left\{\begin{array}{ll}
10r_{m-1}-10&\mbox{if}\quad  r\geq 5,\\
8r_{m-1}-8&\mbox{if}\quad  r=4,\\
6r_{m-1}-6&\mbox{if}\quad  r=3,\\
4r_{m-1}-4&\mbox{if}\quad  r=2.\end{array}\right.
\end{equation}

For this, simply note that $(P_{m-1}+P_m)P_m \geq 2P_{m-1}P_m$ since
$P_{m-1}\leq P_m$ and $P_m$ is nef. Then apply Lemma 2.1. Note that even
if $r_{m-1}=r_1=1$ (hence $P_1=0$) $(P_{m-1}+P_m)P_m \geq 10r_{m-1}-10$
holds.

Finally, if $1\leq i\leq m-2$, then
\begin{equation}\label{eqdosset}
(P_{i}+P_{i+1})P_{m}\geq \left\{\begin{array}{rl}
10r_{i}-5&\mbox{if}\quad  r\geq 5,\\
8r_{i}-4&\mbox{if}\quad  r=4,\\
6r_{i}-3&\mbox{if}\quad  r=3,\\
4r_{i}-2&\mbox{if}\quad  r=2,\end{array}\right.
\end{equation}
which follows immediately from Lemma 5.9, even if $P_i=P_{1}=0$
($r_{1}=1$).

If $\Delta_f=\mbox{deg }{\cal E} =\displaystyle\sum^{\ell}_{i=1}
r_{i}(\mu _{i}-\mu _{i+1})$, call
$\Delta_1=\displaystyle\sum^{m-1}_{i=1} r_{i}(\mu _{i}-\mu _{i+1})$
and $\Delta_2=\Delta_f-\Delta_1$.

\bigskip

Let us prove first (i) and (iii); we can assume then that
$\vert K_{F} \vert$ has
no hyperelliptic subpencil. We get the following inequalities
using (\ref{eqdosquatre}),(\ref{eqdossis}) and
(\ref{eqdosset}):

\begin{eqnarray}
\mbox{If}&&r\geq 5,\nonumber\\
&&K^{3}_{T/B}\geq 10\Delta_1+9\Delta_2-5\mu _{1}-5\mu _{m-1}-7\mu _{m} -4\mu
_{\ell}\geq 9\Delta_f-17\mu _{1}-4\mu _{\ell}.\label{eqdosvuit}\\
\mbox{If}&&r= 4,\nonumber\\
&&K^{3}_{T/B}\geq 8\Delta_1+9\Delta_2-4\mu _{1}-4\mu _{m-1}-9\mu _{m} -4\mu
_{\ell}\geq 8\Delta_1+9\Delta_2-17\mu _{1}-4\mu _{\ell}.\nonumber\\
\mbox{If}&&r=3,\nonumber\\
&&K^{3}_{T/B}\geq 6\Delta_1+9\Delta_2-3\mu _{1}-3\mu _{m-1}-11\mu _{m} -4\mu
_{\ell}\geq 6\Delta_1+9\Delta_2-17\mu _{1}-4\mu _{\ell}.\nonumber
\end{eqnarray}

Note that the bound for $r \geq 5$ also holds for $m=2, \quad P_{m-1}
=P_1=0 \quad \break (
r_1=1), \mbox{or} \quad m=1$.
Using now Remark 2.3 and Lemma 2.1 (i), we have
\begin{eqnarray*}
K^{3}_{T/B}\geq P^{2}_{\ell}(\mu _{1}+2\mu _{\ell})\geq (3p_{g}(F)-7)
(\mu _{1}+2\mu _{\ell})
\end{eqnarray*}
and so (we use  $-17\mu _{1}-4\mu _{\ell}\geq -17(\mu _{1}+2\mu
_{\ell})$; note that  $\mu _{\ell}$ may be zero):
\begin{eqnarray}
&&\mbox{If}\quad r\geq 5\quad \left  (1+\frac{17}{3p_{g}(F)-7}\right)K^{3}_{T/B}\geq
9\Delta_f.\label{eqdosnou}\\
&&\mbox{If}\quad r=4\quad     \left  (1+\frac{17}{3p_{g}(F)-7}\right)K^{3}_{T/B}\geq
8\Delta_1+9\Delta_2.\nonumber\\
&&\mbox{If}\quad r=3\quad     \left
(1+\frac{17}{3p_{g}(F)-7}\right)K^{3}_{T/B}\geq
6\Delta_1+9\Delta_2.\nonumber
\end{eqnarray}

The first inequality proves (i) and holds also when
$m=2, \quad P_{m-1}=P_1=0 \quad (r_1=1), \mbox{or} \quad m=1$.

In order to prove (iii) we can assume from now on that $r=3$ or
$4$, otherwise we have (i) which is stronger than (iii).

We can also assume $m\geq 2$
and, as pointed out, that $\vert P_{m-1}\vert $ is composed with a
pencil.

We divide the argument according to whether the pencil $\vert P_{m-1}\vert $
is irrational or not.

If the pencil is rational we use the same notation as in Lemma
2.1 and Definition 2.2 and set $\widehat D$ for the (possibly singular) general
element of the linear system $\vert \tau _{\ast }P_{i}\vert=\vert
Q_i \vert$ in $F$ (possibly with
base points).

Then using Lemma 2.1 and according to whether the pencil is irrational
or not we have
\begin{eqnarray}
\mbox{for}&&i\leq m-2,\nonumber\\
&&(P_{i}+P_{i+1})(\tau ^{\ast }K_{F})\geq
(4g-4)r_{i}+(2g-2)\quad (\mbox{except if $P_{i}=P_1=0$, $r_{1}=1$)}
\label{eqdosdeu}\\
\mbox{and}&&(P_{i}+P_{i+1})(\tau ^{\ast }K_{F})\geq (4p_{a}(\widehat D)-4-2\widehat
D^{2})r_{i}-(2p_{a}(\widehat D)-2-\widehat D^{2}),\nonumber\\
\mbox{for}&&i=m-1,\nonumber\\
&&(P_{m-1}+P_{m})(\tau ^{\ast
}K_{F})\geq 2P_{m-1}(\tau ^{\ast }K_{F})\geq (4g-4)r_{m-1}\nonumber\\
\mbox{and}&&(P_{m-1}+P_{m})(\tau ^{\ast }K_{F})\geq 2P_{m-1}(\tau ^{\ast }K_{F})
\geq (4p_{a}(\widehat D)-4-2\widehat D^{2})(r_{m-1}-1)\nonumber\\
\mbox{for}&&m\leq i\leq \ell-1,\nonumber\\
&&(P_{i}+P_{i+1})(\tau
^{\ast }K_{F})\geq P_{i}^{2}+P^{2}_{i+1}\nonumber\\
\mbox{and for}&&i=\ell\nonumber\\
&&(P_{\ell}+P^{2}_{\ell+1})(\tau
^{\ast }K_{F})\geq 2P_{\ell}^{2}.\nonumber
\end{eqnarray}

Using Remark 2.3 we know

$$K^{3}_{T/B} \geq \sum^{\ell}_{i=1}(P_{i}+P_{i+1})(\tau ^{\ast
}K_{F})(\mu _{i}-\mu _{i+1})$$
and so we can conclude

\begin{eqnarray}
\mbox{If}&&r=3,4, \, \, p=1\nonumber\\
&&K^{3}_{T/B}\geq (4g-4)\Delta_1+6\Delta_2-11\mu _{m}-3\mu
_{\ell}\geq (4g-4)\Delta_1+6\Delta_2-11\mu _{1}-3\mu _{\ell}.
\label{eqdosonze}
\end{eqnarray}

\vglue.5truecm

Finally:
\begin{eqnarray}
\mbox{if}&&r=3,4, \, \, p=0\nonumber\\
&&K^{3}_{T/B}\geq (4p_{a}(\widehat D)-4-2\widehat
D^{2})\Delta_1+6\Delta_2-(2p_{a}(\widehat D)-2-\widehat D^{2}) (\mu
_{1}-\mu _{m-1})\nonumber\\
&&-(4p_{a}(\widehat D)-4-2\widehat D^{2}) (\mu
_{m-1}-\mu _{m})-11\mu _{m}-3\mu_{\ell}\label{eqdosdotze}\\
&&\geq (2p_{a}(\widehat D)-2-\widehat D^{2})\Delta_1+6\Delta_2-11\mu
_{1}-3\mu _{\ell}\nonumber
\end{eqnarray}
since $\Delta_1\geq (\mu _{1}-\mu _{m-1})+2(\mu _{m-1}-\mu
_{m})$ (this is immediate if $m-1\geq 2$; if $m=2$, then $r_{1}\geq 2$ and
$\Delta_1=r_{1}(\mu _{1}-\mu _{2})\geq 2(\mu _{1}-\mu _{m})$).

Note also that these formulas include the possibility $r_{1}=1$,
$P_{1}=0$.
\vglue.2truecm

Using now that $K^{3}_{T/B} \geq (3p_{g}(F)-7)(\mu_{1}+2\mu_{\ell})$, we get

\begin{eqnarray}
&&\left  (1+\frac{11}{3p_{g}(F)-7} \right)K^{3}_{T/B}\geq
\left\{\begin{array}{ll}
(4g-4)\Delta_1+6\Delta_2&\mbox{if $r=3,4$, $p=1$} \\
(2p_{a}(\widehat D)-2-\widehat D^{2})\Delta_1+6\Delta_2&\mbox{if $r=3,4$,
$p=0.$} \end{array}\right .\nonumber 
\end{eqnarray}

Considering simultaneously this last inequality together with
(\ref{eqdosnou}) and using that $\Delta_2=\Delta_f-\Delta_1$, we get (iii).
All the arguments work similarly
so we only give details of one of them. Assume $P_{m-1}$ is a tetragonal
irrational pencil. Then the last inequality and (14) give

$$\left (1+\frac{17}{3p_g(F)-7}\right )K_{T/B}^3 \geq
8 \Delta_1+9\Delta_2  =8\Delta_f+\Delta_2,$$
$$\left (1+\frac{11}{3p_g(F)-7}\right )K_{T/B}^3 \geq
(4g-4)\Delta_1+6\Delta_2=(4g-4)\Delta_f-(4g-10)\Delta_2.$$

\smallskip

Observe that since $\vert P_{m-1} \vert$ is not an hyperelliptic pencil,
then $g \geq 3$ and so we can get a
lower bound for $\Delta_2$ from the second inequality and we obtain
$$\left [1+\frac{17}{3p_g(F)-7}+\frac{1}{4g-10}
\left (1+\frac{11}{3p_g(F)-7}\right )\right ]K_{T/B}^3 \geq
\left [8 + \frac{4g-4}{4g-10} \right ] \Delta_f$$
which gives $\lambda^1_4(P_{m-1})$. As for the computation of
$\lambda_4^0(P_{m-1})$ or $\lambda_3^0(P_{m-1})$ we only must be
careful when $\delta \leq 6$ since then the corresponding second
inequality does not give a lower bound for $\Delta_2$. If this
case occurs, then just deduce from (14)
$$\left (1+\frac{17}{3p_g(F)-7}\right )K^3_{T/B} \geq 8\Delta_f \quad \quad
\mbox{if the pencil is tetragonal}$$
$$\left (1+\frac{17}{3p_g(F)-7}\right )K^3_{T/B} \geq 6\Delta_f \quad
\mbox{if the pencil is trigonal}$$
which give the special values of $\lambda_r^0$ in (iii).

Of course we must consider all the possibilities for $\vert P_{m-1}
\vert$ being trigonal or tetragonal subpencils of $\vert K_F \vert$
and so we must consider the minimum of all such lower bounds.

\vglue.3truecm

Finally we must prove (iv). Assume $\vert K_{F} \vert$ has hyperelliptic
subpencils. Then it may happen that for some $i \geq m, \quad \varphi_{i}$
is of degree two onto a ruled surface. Also may happen that $r=2$.
Altogether, Remark 2.3, Lemma 2.1 and inequalities (\ref{eqdosvuit}), (\ref{eqdosnou}),
(\ref{eqdosonze}) and (\ref{eqdosdotze}) read

\begin{eqnarray*}
K^{3}_{T/B}\geq P^{2}_{\ell}(\mu _{1}+2\mu _{\ell})\geq (2p_{g}(F)-4)
(\mu _{1}+2\mu _{\ell}).
\end{eqnarray*}
\begin{eqnarray}
\mbox{If}&&r\geq 3\quad \mbox{or} \quad m=1\quad\mbox{or} \quad
m=2 \quad P_{m-1}=P_1=0,\nonumber\\
&&K^{3}_{T/B}\geq 6\Delta_1+6\Delta_2-3\mu _{1}-3\mu _{m-1}-3\mu _{m} -3\mu
_{\ell}\geq 6\Delta_f-9\mu _{1}-3\mu _{\ell}\nonumber
\end{eqnarray}
and so

\begin{eqnarray}
&&\mbox{if}\quad r\geq 3\quad \left  (1+\frac{9}{2p_{g}(F)-4}\right)K^{3}_{T/B}\geq
6\Delta_f.\label{eqhypuno}
\end{eqnarray}

\begin{eqnarray}
\mbox{If}&&r=2,\nonumber\\
&&K^{3}_{T/B}\geq 4\Delta_1+6\Delta_2-2\mu _{1}-2\mu _{m-1}-5\mu _{m} -3\mu
_{\ell}\geq 4\Delta_1+6\Delta_2-9\mu _{1}-3\mu _{\ell}.\nonumber
\end{eqnarray}

And so, if $r=2$, then
\begin{eqnarray}
&&\left
(1+\frac{9}{2p_{g}(F)-4}\right)K^{3}_{T/B}\geq
4\Delta_1+6\Delta_2.\label{eqhypdos}
\end{eqnarray}

If $r=2 \, \, p=0,$
$$K^{3}_{T/B}\geq (2p_{a}(\widehat D)-2-\widehat
D^{2})\Delta_1+4\Delta_2-6\mu _{1}-2\mu _{\ell}.$$

\begin{eqnarray}
\mbox{If}&&r=2\,\, p=1,\nonumber\\
&&K^{3}_{T/B}\geq (4g-4)\Delta_1+4\Delta_2-2\mu _{m}-2\mu _{\ell}
\geq (4g-4)\Delta_1+4\Delta_2-2\mu _{1}-2\mu _{\ell}.\nonumber
\end{eqnarray}

This last inequality needs an extra explanation for the coefficient
of $\mu_m$. If $r=2$, we have an hyperelliptic pencil on ${\widetilde F}$. Let
$D$ be a general irreducible member. For $i \geq m,\, \vert P_i \vert
_{D}$ is a base point free sublinear system of
$\vert K_D \vert$ by adjunction and hence maps $D$ onto ${\proj}^1$.
Hence, $\Sigma_i=\varphi_i({\widetilde F})$ is always a ruled surface.
Moreover, $\varphi_i$ is of degree at least two. Since $\varphi_{m-1}$ factors
through $\Sigma_i$ for $i \geq m$, then $q(\Sigma_i) \geq 1$ (the pencil
is irrational). Hence for $i \geq m$:

$$P_{i}(\tau ^{\ast }K_{F})\geq P_{i}^{2}\geq
(\mbox{deg}\varphi_i)(r_{i}-2+q(\Sigma _{i}))\geq
2r_{i}-2 \, .$$

From here we get

\begin{eqnarray}
&&\left  (1+\frac{2}{2p_{g}(F)-4} \right)K^{3}_{T/B}\geq
(4g-4)\Delta_1+4\Delta_2\qquad \qquad \mbox{if $r=2$,
$p=1$}\label{eqhyptres}\\
\mbox{and}&&\left  (1+\frac{6}{2p_{g}(F)-4} \right)K^{3}_{T/B}\geq
(2p_{a}(\widehat D)-2-\widehat D^{2})\Delta_1+4\Delta_2 \qquad
\mbox{if $r=2$,
$p=0.$}\nonumber
\end{eqnarray}

If $r=2$ (i.e., $\vert P_{m-1} \vert$ is one of the hyperelliptic
subpencils of $\vert K_F \vert$) then we can proceed as in (iii)
using (19) and (20) and inequalities in (iv) follow.
Here the exceptional bounds appear in the rational and irrational cases.
When $p=0$ and $\delta \leq 4$,
(19) gives
$$\left (1+\frac{9}{2p_g(F)-4}\right )K^3_{T/B} \geq 4\Delta_f$$
and so
$$\lambda_1(f) \geq 4\left (1- \frac{9}{2p_g(F)+5}\right ) .$$

When $p=1,\quad g=2$ we have from (20)

$$\lambda_1(f) \geq 4\left (1- \frac{1}{p_g(F)-1}\right ) .$$

If $r \geq 3$ or $\vert P_{m-1} \vert$ is not composed with a pencil
the situation can only be better; indeed, in this case
inequality (18) holds, which is better than inequality (19)
and hence it is better than any inequalities coming from (19)
as those in (iv) are.
\quad  $\Box$

{\bf Remark 2.6.} \quad Theorem 2.4 shows that there is some influence
on the slope
of the existence of certain special maps on the general fibre. This
is precisely what is known to happen partially in the case
of fibred surfaces, where the gonality of a general fibre seems
to play a special role (cf. \cite{SF}, \cite{K8}, \cite{K1}).

{\bf Remark 2.7.} \quad In \cite{B2} using a different aproach
the author finds some lower bounds
for $\lambda_2(f)$ assuming that $T$ is Gorenstein and either
the sheaf ${\cal E}$ is semistable or the canonical image of $F$
lies in few quadrics.

\newpage

\noindent
{\bf  3. The slope of non Albanese fibred threefolds}

\medskip

Consider, for $t\in B$ such that  $F_{t}$ is smooth, the natural diagram
of Albanese maps
$$
\xymatrix{
F_{t}\ar[r]^{\mbox{alb}_{F_{t}}}\ar[d]^{i_{t}}&\mbox{Alb}(F_{t})\ar[d]
_{(i_{t})_{\ast }}\\
T
\ar[r]^{\mbox{alb}_{T}}\ar[d]^{f}&\mbox{Alb}(T)\ar[d]
_{f_{\ast }}\\
B
\ar[r]^{\mbox{alb}_{B}}&\mbox{Alb}(B)}
$$

As $t$ varies, the abelian subvariety $(i_{t})_{\ast}\mbox{Alb}(F_{t})$
remains constant, say $A$, by the rigidity property of abelian varieties.
Let $\alpha =\mbox{dim}A$. From this we get

$$b \leq q(T)=b+\alpha \leq b+g$$

Moreover, from the diagram and the universal property of Albanese
varieties it is clear to see that if $b \geq 1$, then $b=q$ if and only if
$\mbox{alb}_{T}(T)=B$. In this case we say that $f$ is an {\it Albanese
fibration}. We will say that $f$ is a {\it non-Albanese fibration} if
$q(T)>b$ (i.e., if $b=0$ and $q(T) >0$ or if $a=\mbox{dim}alb_{T}(T) \geq 2$).
We want to analyze the influence of this fact in the slope
as in Xiao's result for surfaces (cf. \cite{X1} Corollary 2.1
): if $q(S)>b$ then $\lambda (f) \geq 4$. Our main
tool will be Theorem 2.4 together with an argument of \'etale covers.

First of all we need to control when \'etale covers of curves
are $d$-gonal.

{\bf Lemma 3.1.} \quad
{\it
Let $D$ be a smooth curve and
${\cal L}\in  {\rm Pic}^0(D)$ a  $n$-torsion element.

Let $\alpha :\widetilde{D}\longrightarrow D$
the associated \'etale cover of degree $n$.

Assume  $\widetilde{D}$ has a unique base point free $g^{1}_{d}$;
then

\begin{enumerate}
\item[{\rm (i)}]  $D$ has a $g^{1}_{d}$.
\item[{\rm (ii)}]  $n|d$.
\end{enumerate}
}

\dem Easy.\quad $\Box$

{\bf Lemma 3.2.} \quad
{\it
Let $F$ be a surface of general type and
${\cal L}\in  {\rm Pic}^0(F)$ a  $n$-torsion element. Let $\alpha
:\widetilde{F}\longrightarrow F$ be the associated \'etale cover of degree
$n$. Then, if  $n$ is prime and large enough

\begin{enumerate}
\item[{\rm (i)}] $\widetilde{F}$ has no rational pencil of  $d$-gonal curves
($d=2,3$).
\item[{\rm (ii)}] If $\widetilde{F}$ has an irrational pencil of  $d$-gonal
curves ($d=2,3$) of genus $g$, then so has  $F$ and there exists a base
change
$$
\xymatrix{
\widetilde{F}\ar[d]_{\widetilde{h}}\ar[r]^{\alpha }&F\ar[d]^{h}\\
\widetilde{C}\ar[r]&C
}
$$
such that ${\cal L}=h^{\ast}({\cal M})\in h^{\ast }( {\rm Pic}^0(C))$ and
${\widetilde C} \longrightarrow C$ is induced by ${\cal M}$.
\item[{\rm (iii)}] If $\widetilde{F}$ has a pencil of tetragonal
curves then, either so does $F$ and there exists a base change diagram
as in (ii) (and necessarily the pencil is irrational),
or the pencil $\{\widetilde{D}\}$ is of bielliptic curves
and  $F$ has a pencil $\{D\}$ of bielliptic (hence tetragonal) or
hyperelliptic curves such that $\alpha ^{\ast }(D)=\widetilde{D}$ (and
$2g(\widetilde{D})-2=n(2g(D)-2))$ where
$$
\xymatrix @C=0pt@R=1truecm{
\widetilde{F}\ar[dr]_{\widetilde{h}}&\displaystyle\mathop{\longrightarrow}
^{\alpha }&F\ar[dl]^{h}\\ &C&
}
$$
\end{enumerate}
}

\dem Assume first $\widetilde{F}$ has a base point free pencil
$\widetilde{h}:\widetilde{F}\longrightarrow \widetilde{C}$. Let
$\widetilde{D}$ be a general fibre and let $D=(\alpha
(\widetilde{D}))_{\mbox{red}}$. Clearly $D$ is smooth since
$\widetilde{D}$ moves algebraically.    

If ${\cal L}_{|D}\not= {\cal O}_{D}$, then since  $n$ is prime  ${\cal L}^{\otimes
i}_{|D}\not= {\cal O}_{D}$ for $1\leq i\leq n-1$ and hence $\alpha ^{\ast
}(D)$ is a connected smooth \'etale cover of $D$ containing
$\widetilde{D}$ and so $\alpha ^{\ast }(D)=\widetilde{D}$.

If $D^2>0$, by \cite{MML}, ${\cal L}_{\vert D} \not= {\cal O}_D$ and by
the previous argument $\alpha^{\ast}(D)=\widetilde D$, $0=\widetilde D^2 =
nD^2>0$, which is a contradiction. So necessarily we have $D^2=0$.

If  ${\cal L}_{|D}={\cal O}_{D}$ then  $\alpha ^{\ast
}(D)=\widetilde{D}_{1}+\dots+\widetilde{D}_{n}$ \
($\widetilde{D}_{1}=\widetilde{D}$),
$\widetilde{D}_{i}\widetilde{D}_{j} =0$, $\widetilde{D}_{i}\not=
\widetilde{D}_{j}$ if $i\not= j$, and hence we have a factorization
$$
\xymatrix{
\widetilde{F}\ar[r]^{\alpha }\ar[d]_{\widetilde{h}}&F\ar[d]^{h}\\
\widetilde{C}\ar[r]^{\beta }&C
}
$$
such that  ${\cal L}\in h^{\ast } {\rm Pic}^0(C)$ and $\beta $ is an \'etale
cover. In particular $g(C)\geq 1$.

Now we want to explore the possibility  $\widetilde{D}=\alpha ^{\ast
}(D)$. Since  $n$ is large enough then so is  $g(\widetilde{D})$ and
hence, if  $\widetilde{D}$ has a  $g^{1}_{d}$ ($d=2,3,4$) it is unique
except if $d=4$ and $\widetilde{D}$ is bielliptic. So Lemma
3.1 (ii) says that $\widetilde{D}$ has no  $g^{1}_{d}$
($d=2,3,4$) as long as $n$ does not divide $d$ except when  $\widetilde{D}$  is
bielliptic.

But in this case, let $\sigma $ be the (unique, if
$g(\widetilde{D})\geq 6$) bielliptic involution of $\widetilde{D}$. Let
$\varphi $ be the rank $n$ automorphism of $\widetilde{D}$ induced by
$\alpha $. Since $\sigma $ is unique we must have that $\varphi
\circ\sigma \circ\varphi ^{-1}=\sigma $ and hence there exists an
automorphism $\bar\varphi $ of the elliptic base curve $E$ such that the
following diagram commutes
$$
\xymatrix{
\widetilde{D}\ar[r]^{\varphi }\ar[d]_{2:1}&\widetilde{D}\ar[d]^{2:1}\\
E=\widetilde{D}_{/<\sigma >}\ar[r]^{\bar\varphi }&E}
$$
Hence there is an induced degree two map
$D=\widetilde{D}_{/<\varphi >}\longrightarrow E_{/<\bar \varphi >}=E'$,
where $E'$ is
${\proj}^{1}$ or an elliptic curve according  $\bar \varphi $ has fixed
points or not. Then  $D$ is hyperelliptic or bielliptic (hence
tetragonal).\quad  $\Box$

Then we can consider what is the influence of the irregularity of
$T$ in the slope of $f$. We have an exceptionally good
behaviour:

{\bf Theorem 3.3.} \quad
{\it
Let  $f:T\longrightarrow B$ be a
relatively minimal fibration of a normal,
projective threefold with only canonical singularities onto a smooth curve of
genus $b$. Let $F$ be a general fibre. Assume  $F$ is of general type,
$p_{g}(F)\geq 3$, and that $\chi_f=\chi ({\cal O}_{F})\chi({\cal O}_{B})-
\chi({\cal O}_{T})>0$.

Then, if $q(T)>b$, we have

\begin{enumerate}
\item[{\rm (i)}]  $\lambda _{2}(f)\geq 4$.
\item[{\rm (ii)}] If $F$ has no irrational pencil of  $d$-gonal curves
($d=2,3,4$) then  $\lambda _{2}(f)\geq 9$.
\item[{\rm (iii)}]  If  $F$ has an irrational tetragonal pencil of
genus $g$ and no irrational pencil of trigonal curves nor pencil of
hyperelliptic curves then $\lambda _{2}(f)\geq 9-\frac{3}{4g-9}$.
\item[{\rm (iv)}] If $f$ has an irrational tetragonal pencil of genus
$g_{1}$, an irrational trigonal pencil of genus $g_{2}$ and no
hyperelliptic pencil, then $$\lambda _{2}(f)\geq \min\left\{9-\frac{3}
{4g_{1}-9}, 9-\frac{9}{4g_{2}-7}\right\}.$$
\item[{\rm (v)}] If $F$ has an irrational hyperelliptic pencil of
genus $g$, then  $\lambda _{2}(f)\geq 6-\frac{2}{2g-3}$.
\item[{\rm (vi)}] If  $F$ has a rational hyperelliptic pencil and none
irrational, then $\lambda _{2}(f)\geq 6$.
\item[{\rm (vii)}] If $\lambda _{2}(f)<9$ then either  $F$ has a rational pencil
of hyperelliptic curves or there exists, perhaps up to base change, a
factorization of $f$
$$
\xymatrix{
T\ar@{-->}[r]^{h}\ar[d]^{f}&S\ar[dl]^{g}\\
B&}
$$
where $S$ is a smooth surface fibred over $B$ by curves $C_{t}$ of
genus $g(C_{t})\geq 1$ and $h$ is everywhere defined at the general
fibre $F_{t}$ of $f$, such that
\begin{itemize}
\item for $t\in B$ general the image of $((i_{t})^{\ast }: {\rm Pic}^{
0}(T)\longrightarrow  {\rm Pic}^0(F_{t}))$ lies in $h^{\ast } _{t}
( {\rm Pic}^0(C_{t}))\subseteq  {\rm Pic}^0(F_{t})$.
\item for $s\in S$ general $D_{s}=h^{-1}(s)$ is hyperelliptic, trigonal or tetragonal
(necessarily hyperelliptic or of genus $3$ if $\lambda _{2}(f)<8$).
\end{itemize}
\end{enumerate}
}

\dem Note that  (i) follows from (ii), (iii), (iv), (v) and (vi). Since
$q(T)>b$ then for every $n\gg0$ there exists a $n$-torsion element
${\cal L}\in  {\rm Pic}^0T\setminus f^{\ast}( {\rm Pic}^0(B))$ such that for
$1\leq i\leq n-1$, ${\cal L}^{\otimes i}_{|F}\not={\cal O}_{F}$. Then we can
construct the associated \'etale cover $\alpha
:\widetilde{T}\longrightarrow T$ as in Lemma 1.5 (iii) and
get $\widetilde{f}=f\circ\alpha $ such that $\lambda _{2}(f)=\lambda
_{2}(\widetilde{f})$.

If $\widetilde{F}$ is the fibre of $\widetilde{f}$, we have an \'etale
cover $\alpha _{1}:\widetilde{F}\longrightarrow F$ and hence
$p_{g}(\widetilde{F})\geq \chi({\cal O}_{\widetilde F})=
n\chi({\cal O}_{F})
\geq n$ (note that
$q(\widetilde F)\geq q(F)\geq 1$ since $q(T)>b$).  Since we can do this process
for $n$ as large as is needed, we can take in the bounds of Theorem 2.4
limit when $p_{g}(F)$ goes to infinity.

Assume that $|K_{F}|$ is composed. We have
$p_{g}(\widetilde{F})\geq n$ and either the genus of the
fibre of the pencil or the genus of the base curve increases, except
if $q(F)=1$ and the pencil is elliptic. When
$|K_{F}|$  is composed the pencil can only be rational or elliptic and
the genus of the fibre is at most  $5$ provided $p_{g}\gg0$
(\cite{Be4}, \cite{X5}). So if $n$ is large enough and the pencil is
rational,
$|K_{\widetilde{F}}|$ can not be composed. Since $\lambda
_{2}(f)=\lambda _{2}(\widetilde{f})$ we can assume  $|K_{F}|$ is not
composed with a rational pencil.

Finally, if the pencil is elliptic and $q(F)=1$, note that we can apply
Theorem 2.4 (ii) and get that $\lambda_2(f)\geq 12$ if $g \geq 4$,
$\lambda_2(f)\geq 8$ if $g=3$, $\lambda_2(f)\geq 4$ if $g=2$. So
(ii), (iii), (iv), (v) and (vi) hold. From now on we assume $\vert K_F \vert$
is not composed.

When $F$ has a fibration $h:F \longrightarrow C$ we have an induced
map $\widetilde h=h \circ \alpha_1:F \longrightarrow C$. This map may not
have connected fibres and hence factorizes through an \'etale cover
$\widetilde C \longrightarrow C$. We have two possibilities.

There may exist an unbounded sequence $\{n_i\}_{i \in \Natural} \subseteq
\Natural$ such that for every $i$  \break $\widetilde h_{n_i}$ has connected
fibres over $C$ (hence it is a fibration) or for every $n \geq n_0$,\
$\widetilde h_n$ factorizes through a non trivial \'etale cover
$\widetilde C_n \longrightarrow C$.

In any case we have that, if $g_n$ is the genus of the fibration
$\widetilde F_n \longrightarrow \widetilde C_n$, $g \leq g_n \leq
n(g-1)+1$, the border cases being the two extreme possibilities.

If $C={\proj}^1$ (rational pencil) then ${\widetilde C}_n=C$ for all
$n$ and the sequence $\{ \delta_n \}$ is unbounded. If $F$ has a pencil
of tetragonal curves which are bielliptic, then by Lemma 3.2 (iii)
again it may happen that $\{ g_n \}$ is unbounded. Otherwise by using
Lemma 3.2, we have that $g_n=g$ holds for all $n$.

If the sequence $\{g_n\}$ is bounded, since $\lim_{n \rightarrow \infty}
p_g(\widetilde F_n)=\infty$, we can consider Theorem 2.4 (iii), (iv) and get
the bounds of (iii), (iv), (v) and (vi) (note that $g_n=g$ is the worst
case).

Finally assume $\{g_n\}$ (or $\delta_n$) is unbounded. We must take limit in the bounds
of Theorem 2.4 when $g$ (or $\delta$) and $p_g(F)$ simultaneously (and linearly)
grow. In all the cases, the limit is 9.

If $F$ has no irrational $d$-gonal pencil ($d=2,3,4$), neither has $\widetilde{F}$
by Lemma 5.13. If $F$ has a rational $d-$gonal pencil ($d=2,3,4$), we know
yet that $\lambda_2(f) \geq 9$. So we can assume that if $F$
verifies the hypotheses of Theorem 5.11 (i), then so does
${\widetilde F}$,
and so we get $\lambda_2(f)
\geq 9$ in the limit process. This proves (ii).

In order to prove (vii) note that in the previous arguments we always
have $\lambda _{2}(f)\geq 9$ except when there exists
$h_{t}:F_{t}\longrightarrow C_{t}$ with hyperelliptic, trigonal or
tetragonal fibres (such that $g(C_{t})\geq 1$ when non-hyperelliptic)
and for every ${\cal L}\in  {\rm Pic}^0(T)$ the \'etale cover
$\widetilde{F}_{t}\longrightarrow F_{t}$ given by ${\cal L}_{|F_{t}}$
factorizes through an \'etale cover of $C_{t}$. This says that $\mbox{Im
}((i_{t})^{\ast } {\rm Pic}^0T\longrightarrow  {\rm Pic}^0(F_{t}))$ lies
in the subtorus $h^{\ast }_{t} {\rm Pic}^0C_{t}$.

In order to glue all the maps  $h_{t}$ we can proceed as in Theorem
1.6 (iii).\quad $\Box$

{\bf Corollary 3.4.} \quad
{\it
With the same notations as in Theorem
$3.3$, if $q(T)>b$ then
\begin{enumerate}
\item[{\rm (i)}] If $\lambda _{2}(f)<9$ then $F$ is fibred by hyperelliptic,
trigonal or tetragonal curves.
\item[{\rm (ii)}] If $\lambda _{2}(f)<8$ then $F$ is fibred by
genus $3$ or hyperelliptic curves.
\item[{\rm (iii)}] If $\lambda _{2}(f)<\frac{16}{3}$ then $F$ is fibred by
genus $2$ curves.
\end{enumerate}
}

{\bf Corollary 3.5.} \quad
{\it
With the same notations as in Theorem
$3.3$, if ${\cal E} =f_{\ast }\omega _{T/B}$ has a
quotient of rank one and degree zero, then the same conclusions as in
Theorem $3.3$ hold.

In particular, if  $b=0,1$ and  $F$ is not fibred by  $d$-gonal curves
($d=2,3,4$) and $\lambda _{2}(f)<9$ then ${\cal E} $ is ample.
}

\dem
According to Proposition 1.8 (ii), if ${\cal E}$ has a quotient ${\cal L}$
of rank 1 and degree zero, then it is torsion and so it is trivial after
an \'etale base change $\sigma :{\widetilde B} \longrightarrow B$. Using
now part (i) of Proposition 1.18, the induced new fibration
${\widetilde f}: {\widetilde S} \longrightarrow {\widetilde B}$ verifies
$q({\widetilde T}) > g({\widetilde B})$ and hence Theorem 3.3 applies.
Finally note that both fibrations have the same slope $\lambda_2$ by
Lemma 1.5 (ii).
\quad $\Box$

{\bf Remark 3.6.} \quad
Although in Theorem 2.4 we
only must take care of subpencils of $|K_{F}|$, in the proof of
Theorem 3.3 we must take care of subpencils of
$|K_{\widetilde{F}}|$ for any \'etale cover
$\widetilde{F}\longrightarrow F$, hence corresponding to {\em arbitrary}
pencils in $F$. Hence the hypotheses that appear in the statement of
Theorem 3.3 can not be restricted to subpencils of
$|K_{F}|$.

{\bf Example 3.7.} \quad
We give a family of examples of
fibred threefolds with $F$ fibred by genus two curves and with
slope arbitrarily near to 6. For this, consider a ruled surface $R$
onto a smooth curve $C$ of genus $m$, and let $B$ be a smooth curve
of genus $b$. Let $Y=R \times B$ and consider a suitable double
cover $T \longrightarrow Y$. If the ramification locus is suitably
chosen, a general fibre $F$ of the induced fibration
$f:T\longrightarrow B$ has a genus two fibration. A standard
computation shows that $\lambda_{2}(f)$ is arbitrarily near to 6
provided $m \geq 1$ (in fact equal to 6 if $m=1$). Observe that
by construction $q(T)-b \geq m$ and so $f$ is a non Albanese fibration.
Thus, we can conclude that the bound 6 has certainly some meaning for
fibrations with general fibre fibred by hyperelliptic curves. The
same construction produces examples with arbitrary $g \geq 3$ but
then $\lambda_{2}(f)$ is far from 6.

\vglue.5truecm

\noindent
{ \bf  4. Fibred threefolds with low slope}

\medskip

\indent
In \cite{O}, the following possibilities for fibred threefolds with
fibre of general type and $\lambda_2(f)<4$ are listed.

{\bf Theorem 4.1.}
(Ohno, \cite{O} Main Theorem 2){\bf .} \quad {\it Let $f:T \longrightarrow B$ be a
relatively minimal fibred threefold as in Theorem $2.4$.
Assume $F$ is of general type.
If $K_{T/B}^3 < 4(\chi({\cal O}_B)\chi({\cal O}_F)-\chi({\cal O}_T))$ then
$F$ has one of the following properties:

\begin{itemize}
\item[{\rm (i)}] $F$ carries a linear pencil of curves of genus two.

\item[{\rm (ii)}] $K_F^2 \leq 2 p_g(F)-1$

\item[{\rm (iii)}] $K_F^2=2p_g(F)$, $p_g(F)\geq 3$, $q(F)\leq 2$ and $\vert K_F \vert$
is not composed ($q(F)=2$ only if $p_g(F)=3$).

\item[{\rm (iv)}] $\vert K_F \vert$ is not composed and

\qquad $\cdot$ $K_F^2=8$, $p_g(F)=3$, $q(F)\leq 1$ or

\qquad $\cdot$ $K_F^2=9$, $p_g(F)=4$, $q(F)\leq 1$ or

\qquad $\cdot$ $K_F^2=7$, $p_g(F)=3$, $q(F)\leq 2$

\item[{\rm (v)}] $K_F^2=4$ or $5$, $p_g(F)=2$ and the movable part of $\vert K_F \vert$
is a linear pencil of curves of genus three with only one base point.

\item[{\rm (vi)}] $K_F^2=2$ or $3$ and $p_g(F)=1.$

\item[{\rm (vii)}] $p_g(F)=0.$
\end{itemize}
}

Moreover Ohno gives an example of fibration of type (i).

In the case of fibred surfaces, no examples are known with slope
less than 4 and non-hyperelliptic general fibre if $g=g(F)>>0$.
It seems rather
plausible that if the genus of the fibre is large enough, then
the general fibre must be hyperelliptic.
The analogy in the
case of threefolds is clear: the canonical map $\varphi_{\vert K_F \vert}$
should not be {\it general}. Curiously enough, we can prove this.
If $p_g(F)\leq 2$ the canonical map is
clearly very special. We can prove that if $p_g(F)\geq 8$, $F$ is
fibred by hyperelliptic curves (in fact of genus 2 if $p_g(F) \geq 15)$
and hence $\varphi_{\vert
K_F \vert}$ has at least degree two. In the remaining cases
$3 \leq p_g(F) \leq 7$ we also prove that the canonical map of $F$ has
degree 3 up to some exceptions.

In fact what happens is that only the first case in Ohno's classification
occurs when $p_g(F)$ is large enough. More concretely:

{\bf Theorem 4.2.} \quad
{\it
Let $f:T\longrightarrow B$ be a
relatively minimal fibration of a normal,
projective threefold $T$ with only canonical singularities onto a smooth
curve $B$ of genus $b$. Assume that the general fibre $F$ is of general
type,  $p_{g}(F)\geq 3$ and  $\chi
_{f}=\chi({\cal O}_{F})\chi({\cal O}_{B})-\chi({\cal O}_{T})>0$.

Then, if  $\lambda _{2}(f)<4$, we have:
\begin{enumerate}
\item[{\rm (i)}]  $q(T)=b$
\item[{\rm (ii)}]  ${\cal E} =f_{\ast }\omega_{T/B}$  has no invertible rank
zero quotient sheaf (in particular,  ${\cal E} $ is ample provided
$b\leq 1$).
\item[{\rm (iii)}] If $p_g(F) \geq 15$ then $F$ has a rational pencil of curves
of genus $2$.

\item[{\rm (iv)}] If $p_g(F) \leq 14$ then one of the following holds
\begin{itemize}
\item[{\rm (a)}]  $F$ has a rational pencil of hyperelliptic curves
\item[{\rm (b)}]  $F$ has a rational pencil of trigonal curves, $q(F)=0$ and
\begin{itemize}
\item[.] either the canonical map of $F$ is of degree $3$ and either
$p_g(F)=3, \quad 3 \leq K^2_F \leq 8$ or $p_g(F)=4,5 \quad
3p_g(F)-6 \leq K^2_F \leq 9$
\item[.] or the canonical map of $F$ is birational and

$$p_g(F)=4 \quad 5 \leq K^2_F \leq 9$$
$$5 \leq p_g(F) \leq 7 \quad 3p_g(F)-7 \leq K^2_F \leq 2p_g(F).$$
\end{itemize}

\item[{\rm (c)}]  $F$ is the quintic surfce in $\proj^{3}$ ( that is, $F$ is
canonical, $p_{g}(F)=4$, $q(F)=0$,  $K^{2}_{F}=5$).
\end{itemize}
\end{enumerate}
}

{\bf Remark 4.3.} \quad
It is doubtful that the cases of fibre $F$ canonical in (iv)(b) and in
(iv)(c) occur. In
\cite{B2} the author proves that then
$K^{3}_{T/B}\geq 4\chi_f$, provided $T$ is Gorenstein, so any
example should necessarily have $T$ non Gorenstein.

\dem

The first two statements follow from Theorem 3.3 and Corollary 3.5.
Following the list of Ohno in Theorem 4.1, if $\lambda_2(f)
<4$ and $p_g(F)\geq 3$, $|K_F|$ is not composed.

We follow the notations of the proof of Theorem 2.4.

Assume $F$ has a rational hyperelliptic pencil. We must prove that
the pencil is of genus 2 provided $p_g(F) \geq 15$. Put
$\delta=K_F{\widehat D}$. If the (geometric) genus of ${\widehat D}$ is
not 2 then we observe that $\delta \geq 4$ (if ${\widehat D}^2=0$ then
$\delta=2g-2$; if ${\widehat D}^2>0$ Hodge index theorem gives
$\delta^2 \geq K^2_F \geq 2p_g(F)-4 \geq 26$). Formula (17) reads
for $r=2, p=0$

$$K^3_{T/B} \geq 2\delta \Delta_1 +4\Delta_2 -\delta(\mu_1-\mu_{m-1})
-2\delta(\mu_{m-1}-\mu_{m})-6\mu_m-2\mu_{\ell}$$
$$\geq 2\delta \Delta_1+4\Delta_2-2 \delta \mu_1$$
and using that $K^3_{T/B} \geq (2p_g(F)-4)(\mu_1+2\mu_{\ell})$ we get

$$(1+ \frac{2 \delta}{2p_g(F)-4})K^3_{T/B} \geq 2\delta \Delta_1
+4\Delta_2$$
which together with (19) gives that $K^3_{T/B} \geq 4\Delta_f$ provided
$p_g(F) \geq 15$.

Assume $F$ has no rational hyperelliptic pencil. According to Remark 2.3 we
must check when the coefficient of $(\mu_i - \mu_{i+1})$ is greater
than or
equal to $4r_i$.

Take $i$ such that $m \leq i \leq \ell -1$. If $a_i=2$, then $P_i^2\geq
2r_i-2$ if the image is ruled (since $F$ has no hyperelliptic
rational pencil, deg$\varphi_i(F)\geq r_i-1$) or $P_i^2 \geq 4r_i-8$
otherwise. In any case $P_i^2 \geq 2r_i-2$ ($r_i \geq 3$ since
$\vert P_i \vert$ is not composed), and hence
$P_i^2+P_iP_{i+1}+P_{i+1}^2\geq 3P_i^2\geq 6r_i-6\geq 4r_i$.

If $a_i=3$, then $a_{i+1}=1$ or $3$ and hence $P_{i+1}^2\geq
3r_{i+1}-7\geq 3r_i-4$. If $P_i\not= P_{i+1}$, then $P_i^2+P_iP_{i+1}+
P_{i+1}^2\geq 2P_i^2+1+P_{i+1}^2\geq 9r_i-15\geq 4r_i$.
If $P_i=P_{i+1}$, then $P_i^2+P_iP_{i+1}+P_{i+1}^2=3P_{i+1}^2\geq
9r_i-12 \geq 4r_i$.

If $a_i\geq 4$, $P_i^2+P_iP_{i+1}+P_{i+1}^2\geq 3P_i^2\geq 12r_i-24\geq
4r_i$. Finally if $a_i=1$, then $a_{i+1}=1$ and hence, by the same
argument as in (\ref{eqdosquatre}), $P_i^2+P_iP_{i+1}+P_{i+1}^2\geq 9r_i-17\geq 4r_i$
(since $r_i\geq 4$ if $\varphi_i$ is birational).

Let $i=\ell$. As pointed out in Remark 2.3,
we can set $P_{\ell +1}=\tau^{\ast}K_F$. Hence, if $P_{\ell}=
P_{\ell +1}=\tau^{\ast}K_F$, we have
$P_{\ell}^2+P_{\ell}P_{\ell +1}+P_{\ell +1}^2=3P_{\ell}^2$ but if
$P_{\ell +1} \neq P_{\ell}$, we have $P_{\ell}^2+P_{\ell}P_{\ell +1}+
P_{\ell +1}^2 \geq K^2_F+2P_{\ell}^2+1 \geq 3P_{\ell}^2+2$. Keeping
this in mind we obtain that
$P_{\ell}^2+P_{\ell}P_{\ell +1}+P_{\ell +1}^2
\geq 4 r_\ell$ except when $r_\ell = p_g(F)=4$, $P_\ell^2=K_F^2=3p_g
-7=5$ and $F$ is canonical or $r_\ell=p_g(F)=3$, $P_\ell^2=3p_g-6=3$,
$K^2_F=3p_g(F)-6 \quad \mbox{or} \quad 3p_g(F)-5$ and the canonical map is of
degree three. In both cases we necessarily have $q(F)=0$ (see
\cite{B2} for the canonical case and \cite{Z3} for the degree 3 case).

Take $i$ such that $1 \leq i \leq m-1$. If $r_i=1$
(then $i=1$ and $P_1=0$), we have
$(P_1+P_2)P_m\geq 4r_1=4$ except when $r_2=2$ and $P_m$ induces a
$g_3^1$ in the fibre of the rational pencil $\vert P_2\vert$. Assume $r_i \geq 2$.
If $a_m=2$, then $P_iP_m \geq 2r_i$; for this we must look at the
proof of Lemma 2.1 (ii). Assume
$2r_i-1 \geq P_i P_m \geq \alpha_2(\alpha_1 ad) \geq
(\alpha_1 ad)(r_i-1)$;
we have that $\vert P_m \vert_{\vert D}=g_2^1$, hence
$a=2$ and necessarily $\alpha_1=d=1$; if $\alpha_2=r_i-1$ the pencil
would be rational (since $\alpha_1 =1$) which is impossible by our
assumptions; hence $\alpha_2 \geq r_i$ which
is again impossible.
Then $(P_i+P_{i+1})P_m \geq 2 P_iP_m \geq 4r_i$.
If $a_m \geq 3$, then by Lemma 2.1 \quad $P_iP_m \!\! \geq \!\! 3(r_i-1)\geq 2r_i$ except if $r_i=2$,
$a_m=3$. In this exceptional case, if $P_{i+1}\not= P_m$, then
$(P_i+P_{i+1})P_m \geq (3r_i-3)+(3 r_{i+1}-3)\geq 6r_i-3=9>8=4r_i$; if
$P_{i+1}=P_m$, then $(P_i+P_m)P_m\geq 8=4r_i$ except if
$4 \geq P_m^2\geq 3r_m-6$, i.e., $r_1=2$, $r_2=3$, $m=2$, $a_m=3$
(which again produces a rational trigonal pencil in $F$). Finally if
$a_m=1$, then $P_iP_m \geq 4 r_i-4$ (Lemma 2.1) and hence
$P_iP_m\geq 2 r_i$ and
$(P_i+P_{i+1})P_m\geq 4 r_i$.

So we can conclude than either $p_g(F)=4$, $q(F)=0$, $K_F^2=5$ or
$F$ has a rational trigonal pencil.

Note that in the discussion above, when $F$ has a rational trigonal
pencil, $\vert P_m\vert $ induces a degree 3 map. Hence the canonical map of $F$ can
only be of degree 1 or 3. In any case $K_F^2 \geq 3p_g(F)-7$. Hence,
applying Theorem 4.1 we have $3p_g(F)-7\leq K_F^2 \leq 2p_g(F)$
(if $p_g(F)\geq 5$) and hence $p_g(F)\leq 7$, $K^{2}_{F} \leq
2p_g(F) \leq 14$.

Finally we must prove that $q(F)=0$. If $q(F)=1$, then $K_F^2 \geq 3 p_g(F)+7q(F)-7=3p_g(F)$
(cf. \cite{K4}) which is impossible. Assume $q(F)\geq 2$. If $\vert K_F \vert$
is birational, we have $3p_g(F)-4 \leq K_F^2 \leq 2p_g(F)-1$
(if $q(F) \geq 1$, we have $K^2_F \geq 3p_g(F)+q(F)-7$
but if equality holds then $q(F) \geq 3$ (cf., e.g., \cite{B2}))
which is impossible. If $\vert K_F \vert$ induces a
map of degree 3, we have $3p_g(F)-3 \leq K_F^2 \leq 2p_g(F)$ (cf.
\cite{Z3} and Theorem 4.1) or $p_g(F)=3$, $q(F)=2$, $K_F^2=7$;
so in any case
we get $p_g(F)=3$, $K_F^2\geq 6$. Following the above discussion
the only possibilities for the Harder-Narasimhan filtration of
${\cal E}$ are $r_1=2$, $r_2=3$ or $r_1=1$, $r_2=2$, $r_3=3$.
The first one gives
$$K_{T/B}^3 \geq (P_1+P_2)P_2 (\mu_1-\mu_2)+3
P_2^2\mu_2 \geq 9 (\mu_1 - \mu_2)+18\mu_2 \geq 8 (\mu_1 - \mu_2)+12
\mu_2=4 \Delta_f.$$

The last one gives
$$\begin{array}{rl}
K_{T/B}^3 &\geq (P_1+P_2)P_3
(\mu_1-\mu_2)+(P_2+P_3)P_3(\mu_2 - \mu_3)+3P_3^2\mu_3 \\
&\geq
3 (\mu_1 - \mu_2)+9(\mu_2 - \mu_3)+18 \mu_3 \geq 4 (\mu_1 - \mu_2)+
8(\mu_2- \mu_3)+12 \mu_3=4 \Delta_f
\end{array}$$
if $\mu_2 - \mu_3 \geq \mu_1 - \mu_2$; otherwise consider
$$\begin{array}{rl}
K_{T/B}^3 &\geq (P_1+P_3)P_3(\mu_1 - \mu_3)+
3P_3^2\mu_3 \\
&\geq 6 (\mu_1 - \mu_3)+18\mu_3 \geq 4 (\mu_1 - \mu_2)+
8(\mu_2- \mu_3)+12 \mu_3 = 4 \Delta_f.
\end{array}$$

So we necessarily have $q(F)=0$.

As for the restrictions for $(p_g(F),K^{2}_{F})$ when the canonical
map is of degree 3, we refer to \cite{Po}, \cite{Z3}, \cite{K3}.
\quad $\Box$

\end{document}